\documentclass[10pt,a4paper]{amsart}%
\usepackage{amsfonts,amsmath,amssymb}
\usepackage{hyperref}
\usepackage[all]{xy}

\newtheorem*{theorem*}{Theorem}
\newtheorem{lemma}{Lemma}[subsection]
\newtheorem{proposition}[lemma]{Proposition}
\newtheorem{remark}[lemma]{Remark}

\newtheorem{theorem}[lemma]{Theorem}
\newtheorem{definition}[lemma]{Definition}
\newtheorem{notation}[lemma]{Notation}

\newtheorem{property}[lemma]{Property}
\newtheorem{corollary}[lemma]{Corollary}

\newtheorem{conjecture}{Conjecture}
\newtheorem*{conjecture*}{Conjecture}

\baselineskip 24pt \textwidth 15 cm  \theoremstyle{plain}
\newcommand{\tr}{\operatorname{Tr}}

\newcommand{\cc}{\mathbb{C}}

\newcommand{\eps}{\varepsilon}

\newcommand{\re}{\operatorname{Re}}

\newcommand{\Z}{{\mathbb Z}}

\newcommand{\R}{{\mathbb R}}
\newcommand{\C}{{\mathbb C}}

\newcommand{\E}{{\mathcal E}}

\newcommand{\Sc}{{\mathcal S}}
\newcommand{\G}{{\mathcal G}}

\newcommand{\Fre}{{Fr\'{e}chet \,}}
\newcommand{\et}{{\`{e}tale }}

\newcommand{\Fou}{{\mathcal{F}}}

\newcommand{\cD}{{\mathcal{D}}}
\newcommand{\g}{{\mathfrak{g}}}
\newcommand{\h}{{\mathfrak{h}}}
\newcommand{\z}{{\mathfrak{z}}}

\newcommand{\Supp}{\mathrm{Supp}}

\newcommand{\gd}{\g^{\sigma}}
\newcommand{\oF}{{\overline{F}}}
\newcommand{\cO}{{\mathcal{O}}}
\newcommand{\Sym}{\operatorname{Sym}}

\begin{document}

\author{Avraham Aizenbud}
\address{Avraham Aizenbud and Dmitry Gourevitch, Faculty of Mathematics
and Computer Science, The Weizmann Institute of Science POB 26,
Rehovot 76100, ISRAEL.} \email{aizenr@yahoo.com}
\author{Dmitry Gourevitch} \email{guredim@yahoo.com}

\title[Generalized Harish-Chandra descent]{Generalized Harish-Chandra descent and applications to Gelfand pairs}
\keywords{Multiplicity one, Gelfand pairs, symmetric pairs, Luna
slice theorem, invariant distributions, Harish-Chandra descent. \\
\indent MSC Classes: 20C99, 20G05, 20G25, 22E45, 46F10, 14L24,
14L30}
%
%
%
%
%
%
%
%
%
%
\maketitle
$\quad \quad \quad$ with appendix A by Avraham Aizenbud, Dmitry
Gourevitch and Eitan Sayag

\begin{abstract}
In the first part of the paper we generalize a descent technique
due to Harish-Chandra to the case of a reductive group acting on a
smooth affine variety both defined over arbitrary local field $F$
of characteristic zero. Our main tool is Luna slice theorem.

In the second part of the paper we apply this technique to
symmetric pairs. In particular we prove that the pair
$(GL_n(\C),GL_n(\R))$ is a Gelfand pair. We also prove that any
conjugation invariant distribution on $GL_n(F)$ is invariant with
respect to transposition. For non-archimedean $F$ the later is a
classical theorem of Gelfand and Kazhdan.

We use the techniques developed here in our subsequent work
\cite{AG3} where we prove an archimedean analog of the theorem on
uniqueness of linear periods by H. Jacquet and S. Rallis.
\end{abstract}

\setcounter{tocdepth}{2}
 \tableofcontents

\section{Introduction}
Harish-Chandra developed a technique based on Jordan decomposition
that allows to reduce certain statements on conjugation invariant
distributions on a reductive group to the set of unipotent
elements, provided that the statement is known for certain
subgroups (see e.g. \cite{HCh}).

In this paper we generalize part of this technique to the setting
of a reductive group acting on a smooth affine algebraic variety,
using Luna slice theorem. Our technique is oriented towards
proving Gelfand property for pairs of reductive groups.

Our approach is uniform for all local fields of characteristic
zero - both archimedean and non-archimedean.

\subsection{Main results}$ $\\
The core of this paper is Theorem \ref{Gen_HC}:
\begin{theorem*}
Let a reductive group $G$ act on a smooth affine variety $X$, both
defined over a local field $F$ of characteristic zero. Let $\chi$
be a character of $G(F)$.

Suppose that for any $x \in X(F)$ with closed orbit there are no
non-zero distributions on the normal space to the orbit $G(F)x$ at
$x$ which are equivariant with respect to the stabilizer of $x$
with the character $\chi$.

Then there are no non-zero $(G(F),\chi)$-equivariant distributions
on $X(F)$.
\end{theorem*}

Using this theorem we obtain its stronger version (Corollary
\ref{Strong_HC_Cor}). This stronger version is based on an
inductive argument which shows that it is enough to prove that
there are no non-zero equivariant distributions on the normal
space to the orbit $G(F)x$ at $x$ under the assumption that all
such distributions are supported in a certain closed subset which
is an analog of the cone of nilpotent elements.

Then we apply this stronger version to problems of the following
type. Let a reductive group $G$ acts on a smooth affine variety
$X$, and $\tau$ be an involution of $X$ which normalizes the
action of $G$. We want to check whether any $G(F)$-invariant
distribution on $X(F)$ is also $\tau$-invariant. Evidently, there
is the following
necessary condition on $\tau$:\\
(*) Any closed orbit in $X(F)$ is $\tau$-invariant.\\
In some cases this condition is also sufficient. In these cases we
call the action of $G$ on $X$ \emph{tame}.

The property of being tame is weaker than the property called
"density" in \cite{RR}. However, it is sufficient for the purpose
of proving Gelfand property for pairs of reductive groups.

In section \ref{SecTame} we give criteria for tameness of actions.
In particular, we have introduced the notion of "special" action.
This notion can be used in order to show that certain actions are
tame (see Theorem \ref{Invol_HC} and Proposition
\ref{SpecWeakReg}). Also, in many cases one can verify that an
action is special using purely algebraic - geometric means.

Then we restrict our attention to the case of symmetric pairs.
There we introduce a notion of regular symmetric pair (see
Definition \ref{DefReg}), which also helps to prove Gelfand
property. Namely, we prove Theorem \ref{GoodHerRegGK}.

\begin{theorem*}
Let $G$ be a reductive group defined over a local field $F$ and
$\theta$ be an involution of $G$. Let $H:=G^{\theta}$ and let
$\sigma$ be the anti-involution defined by
$\sigma(g):=\theta(g^{-1})$. Consider the symmetric pair $(G,H)$.

Suppose that all its "descendants" (including itself, see
Definition \ref{descendant}) are regular. Suppose also that any
closed $H(F)$-double coset in $G(F)$ is $\sigma$- invariant.

 Then every
$H(F)$ double invariant distribution on $G(F)$ is
$\sigma$-invariant. In particular, the pair $(G,H)$ is a Gelfand
pair (see section \ref{Gel}).
\end{theorem*}

Also, we formulate an algebraic-geometric criterion for regularity
of a pair (Proposition \ref{SpecCrit}).

Using our technique we prove (in section \ref{Sec2RegPairs}) that
the pair $(G(E),G(F))$ is \emph{tame} for any reductive group $G$
over $F$ and a quadratic field extension $E/F$. This means that
the two-sided action of $G(F)\times G(F)$ on $G(E)$ is tame. This
implies that the pair $(GL_n(E),GL_n(F))$ is a Gelfand pair. In
the non-archimedean case this was proven in \cite{Fli}.

Also we prove that the adjoint action of a reductive group on
itself is tame. This is a generalization of a classical theorem by
Gelfand and Kazhdan, see \cite{GK}.

In our subsequent work \cite{AG3} we use the results of this paper
to prove that the pair $(GL_{n+k},GL_n \times GL_k)$ is a Gelfand
pair by proving that it is regular. In the non-archimedean case
this was proven in \cite{JR} and our proof follows their lines.

In general, we conjecture that any symmetric pair is regular. This
would imply van Dijk conjecture:

\begin{conjecture*}[van Dijk]
Any symmetric pair $(G,H)$ over $\C$ such that $G/H$ is connected
is a Gelfand pair.
\end{conjecture*}

\subsection{Related works on this topic}
$ $\\This paper was inspired by the paper \cite{JR} by Jacquet and
Rallis where they prove that the pair $(GL_{n+k}(F),GL_n(F) \times
GL_k(F))$ is a Gelfand pair for non-archimedean local field $F$ of
characteristic zero. Our aim was to see to what extent their
techniques generalize.

Another generalization of Harish-Chandra descent using Luna slice
theorem has been done in the non-archimedean case in \cite{RR}. In
that paper Rader and Rallis investigated spherical characters of
$H$-distinguished representations of $G$ for symmetric pairs
$(G,H)$ and checked the validity of what they call "density
principle" for rank one symmetric pairs. They found out that
usually it holds, but also found counterexamples.

In \cite{vD}, van-Dijk investigated rank one symmetric pairs in
the archimedean case and gave the full answer to the question
which of them are Gelfand pairs. In \cite{Bos-vD}, van-Dijk and
Bosman studied the non-archimedean case and gave the answer for the
same question for most rank one symmetric pairs. We hope that the second part of our
paper will enhance the understanding of this question for
symmetric pairs of higher rank.

\subsection{Structure of the paper} $ $\\

In section \ref{Prel} we introduce notation that allows us to
speak uniformly about spaces of points of smooth algebraic
varieties over archimedean and non-archimedean local fields, and
equivariant distributions on those spaces.

In subsection \ref{PrelLoc} we formulate a version of Luna slice
theorem for points over local fields (Theorem \ref{LocLuna}). In
subsection \ref{PrelDist} we formulate theorems on equivariant
distributions and equivariant Schwartz distributions.

In section \ref{SecDescent} we formulate and prove the generalized
Harish-Chandra descent theorem and its stronger version.

Section \ref{DistVerSch} is relevant only
to the archimedean case. In that section we prove that in cases that we
consider if there are no equivariant Schwartz distributions then
there are no equivariant distributions at all.
Schwartz
distributions are discussed in Appendix \ref{AppSubFrob}.

In section \ref{SecFour} we formulate homogeneity theorem that
helps us to check the conditions of the generalized Harish-Chandra
descent theorem. In the non-archimedean case this theorem had been
proved earlier (see e.g. \cite{JR}, \cite{RS2} or \cite{AGRS}). We
provide the proof for the archimedean case in Appendix
\ref{AppRealHom}.

In section \ref{SecTame} we introduce the notion of tame actions
and provide tameness criteria.

In section \ref{SecSymPairs} we apply our tools to symmetric
pairs. In subsection \ref{SecTamePairs} we provide criteria for
tameness of a symmetric pair. In subsection \ref{SecRegPairs} we
introduce the notion of regular symmetric pair and prove Theorem
\ref{GoodHerRegGK} that we mentioned above. In subsection
\ref{conj} we discuss conjectures about regularity and Gelfand
property of symmetric pairs. In subsection \ref{Sec2RegPairs} we
prove that certain symmetric pairs are tame.

In section \ref{Gel} we give preliminaries on Gelfand pairs an
their connections to invariant distributions. We also prove that
the pair $(GL_n(E), GL_n(F))$ is Gelfand pair for any quadratic
extension $E/F$.

In Appendix \ref{SecRedLocPrin} we formulate and prove a version
of Bernstein's localization principle (Theorem \ref{LocPrin}).
This is relevant only for archimedean $F$ since for $l$-spaces
a more general version of this principle had been proven in
\cite{Ber}.  This appendix is used in section \ref{DistVerSch}.

In \cite{AGS2} we formulated localization principle in the setting
of differential geometry. Currently we do not have a proof of this
principle in such general setting. In Appendix \ref{SecRedLocPrin}
we present a proof in the case of a reductive group $G$ acting on
a smooth affine variety $X$. This generality is wide enough for
all applications we had up to now, including the one in
\cite{AGS2}.

We start Appendix \ref{AppLocField} from discussing different
versions of the inverse function theorem for local fields. Then we
prove a version of Luna slice theorem for points over local fields
(Theorem \ref{LocLuna}). For archimedean $F$ it was done by Luna
himself in \cite{Lun2}.

Appendices \ref{AppSubFrob} and \ref{AppRealHom} are relevant only
to the archimedean case.

In Appendix \ref{AppSubFrob} we discuss Schwartz distributions on
Nash manifolds. We prove for them Frobenius reciprocity and
construct a pullback of a Schwartz distribution under Nash
submersion. Also we prove that $K$ invariant distributions which
are (Nashly) compactly supported modulo $K$ are Schwartz distributions.

In Appendix \ref{AppRealHom} we prove the archimedean version of
the homogeneity theorem discussed in section \ref{SecFour}.

In Appendix \ref{Diag} we present a diagram that illustrates the
interrelations of various properties of symmetric pairs.

\subsection{Acknowledgements}
We would like to thank our teacher \textbf{Joseph Bernstein} for
our mathematical education.

We also thank \textbf{Vladimir Berkovich}, \textbf{Joseph
Bernstein}, \textbf{Gerrit van Dijk}, \textbf{Stephen Gelbart},
\textbf{Maria Gorelik}, \textbf{David Kazhdan}, \textbf{Erez
Lapid}, \textbf{Shifra Reif}, \textbf{Eitan Sayag}, \textbf{David
Soudry}, \textbf{Yakov Varshavsky} and \textbf{Oksana Yakimova}
for fruitful discussions, and \textbf{Sun Binyong} for useful
remarks.

Finally we thank \textbf{Anna Gourevitch} for the graphical design
of Appendix \ref{Diag}.

\section{Preliminaries and notations} \label{Prel}

\begin{itemize}
\item
From now and till the end of the paper we fix a local field $F$ of
characteristic zero. All the algebraic varieties and algebraic
groups that we will consider will be defined over $F$.
\item For a group $G$ acting on a set $X$ and an element $x \in X$
we denote by $G_x$ the stabilizer of $x$.
\item By a reductive group we mean an algebraic reductive group.

\end{itemize}

We treat an algebraic variety $X$ defined over $F$ as algebraic
variety over $\oF$ together with action of the Galois group
$Gal(\oF,F)$. On $X$ we will consider only the Zariski topology.
On $X(F)$ we consider only the analytic (Hausdorff) topology. We
treat finite dimensional linear spaces defined over $F$ as
algebraic varieties.

Usually we will use letters $X, Y, Z, \Delta$ to denote algebraic
varieties and letters $G, H$ to denote algebraic groups. We will
usually use letters $V, W,U,K,M,N,C,O,S,T$ to denote analytic
spaces and in particular $F$ points of algebraic varieties and the
letter $K$ to denote analytic groups. Also we will use letters $L,
V, W$ to denote vector spaces of all kinds.

\setcounter{lemma}{0}
\begin{definition}
Let an algebraic group $G$ act on an algebraic variety $X$. A pair
consisting of an algebraic variety $Y$ and a $G$-invariant
morphism $\pi:X \to Y$ is called \textbf{the quotient of $X$ by
the action of $G$} if for any pair $(\pi ', Y')$, there exists a
unique morphism $\phi : Y \to Y'$ such that $\pi ' = \phi \circ
\pi$. Clearly, if such pair exists it is unique up to canonical
isomorphism. We will denote it by $(\pi _X, X/G)$.
\end{definition}

\begin{theorem} \label{Quotient}
Let a reductive group $G$ act on an affine variety $X$. Then the
quotient $X/G$ exists, and every fiber of the quotient map $\pi_X$
contains a unique closed orbit.
\end{theorem}
\begin{proof}
In \cite{Dre} it is proven that the variety $Spec\cO(X)^G$
satisfies the universal condition of $X/G$. Clearly, this variety
is defined over $F$ and hence we can take $X/G:=Spec\cO(X)^G$.
\end{proof}

\subsection{Preliminaries on algebraic geometry over local fields}
\label{PrelLoc}

\subsubsection{Analytic manifolds}
$ $\\
In this paper we will consider distributions over $l$-spaces,
smooth manifolds and Nash manifolds. $l$-spaces are locally
compact totally disconnected topological spaces and Nash manifolds
are semi-algebraic smooth manifolds.

For basic preliminaries on $l$-spaces and distributions over them
we refer the reader to \cite{BZ}, section 1.

For preliminaries on Nash manifolds and Schwartz functions and
distributions over them see Appendix \ref{AppSubFrob} and
\cite{AG1}. In this paper we will consider only separated Nash
manifolds.

We will now give notations which will allow a uniform exposition
of archimedean and non-archimedean cases.

We will use the notion of analytic manifold over a local field
(see e.g. \cite{Ser}, Part II, Chapter III). When we say
"\textbf{analytic manifold}" we mean analytic manifold over some
local field. Note that an analytic manifold over a non-archimedean
field is in particular an $l$-space and analytic manifold over an
archimedean field is in particular a smooth manifold.

\begin{definition}
A \textbf{B-analytic manifold} is either an analytic manifold over
a non-archimedean local field, or a Nash manifold.
\end{definition}

\begin{remark} If $X$ is a smooth algebraic variety,
then $X(F)$ is a B-analytic manifold and $(T_xX)(F) = T_x(X(F)).$
\end{remark}

\begin{notation}
Let $M$ be an analytic manifold and $S$ be an analytic
submanifold. We denote by $N_S^M:=(T_M|_Y)/T_S $ the normal bundle
to $S$ in $M$. The conormal bundle is defined by
$CN_S^M:=(N_S^M)^*$. Denote by $Sym^k(CN_S^M)$ the k-th symmetric
power of the conormal bundle. For a point $y\in S$ we denote by
$N_{S,y}^M$ the normal space to $S$ in $M$ at the point $y$ and by
$CN_{S,y}^M$ the conormal space.
\end{notation}

\subsubsection{$G$-orbits on $X$ and $G(F)$-orbits on $X(F)$} 

\begin{lemma} \label{OrbitIsOpen}
Let $G$ be an algebraic group. Let $H \subset G$ be a closed
subgroup. Then $G(F)/H(F)$ is open and closed in $(G/H)(F)$.
\end{lemma}
For proof see Appendix \ref{AppSub}.

\begin{corollary}
Let an algebraic group $G$ act on an algebraic variety $X$. Let $x
\in X(F)$. Then $$N_{Gx,x}^{X}(F) \cong N_{G(F)x,x}^{X(F)}.$$
\end{corollary}

\begin{proposition} \label{LocClosedOrbit}
Let an algebraic group $G$ act on an algebraic variety $X$.
Suppose that $S \subset X(F)$ is non-empty closed $G(F)$-invariant
subset. Then $S$ contains a closed orbit.
\end{proposition}
\begin{proof}
The proof is by Noetherian induction on $X$. Choose $x \in S$.
Consider $Z:=\overline{Gx} -Gx$.

If $Z(F)\cap S$ is empty then $Gx(F) \cap S$ is closed and hence
$G(F)x \cap S$ is closed by Lemma \ref{OrbitIsOpen}. Therefore
$G(F)x$ is closed.

If $Z(F)\cap S$ is non-empty then $Z(F) \cap S$ contains a closed
orbit by the induction assumption.
\end{proof}

\begin{corollary} \label{OpenClosedAll}
Let an algebraic group $G$ act on an algebraic variety $X$. Let
$U$ be an open $G(F)$-invariant subset of $X(F)$. Suppose that it
includes all closed $G(F)$-orbits. Then $U=X(F)$.
\end{corollary}

\begin{theorem} \label{LocZarClosed}
Let a reductive group $G$ act on an affine variety $X$. Let $x
\in X(F)$. Then the following are equivalent:\\
(i) $G(F)x \subset X(F)$ is closed (in the analytic topology).\\
(ii) $Gx\subset X$ is closed (in the Zariski topology).
\end{theorem}
For proof see \cite{RR}, section 2 fact A, pages 108-109.

\begin{definition}
Let a reductive group $G$ act on an affine variety $X$. We call an
element $x \in X$ \textbf{$G$-semisimple} if its orbit $Gx$ is
closed. In particular, in the case of $G$ acting on itself by the
adjoint action, the notion of $G$-semisimple element coincides
with the usual notion of semisimple element.
\end{definition}

\begin{notation}
Let $V$ be an algebraic finite dimensional representation over $F$
of a reductive group $G$. We denote
$$Q(V):=(V/V^G)(F).$$ Since $G$ is reductive, there is a
canonical embedding $Q(V) \hookrightarrow V(F)$. Let $\pi : V(F)
\to (V/G)(F)$ be the standard projection. We denote $$\Gamma(V):=
\pi^{-1}(\pi(0)).$$ Note that $\Gamma(V) \subset Q(V)$. We denote
also $$R(V):= Q(V) - \Gamma(V).$$
\end{notation}

\begin{notation}
Let a reductive group $G$ act on an affine variety $X$. Let an
element $x \in X(F)$ be $G$-semisimple. We denote $$S_x := \{y \in
X(F) \, | \, \overline{G(F)}y \ni x\}.$$
\end{notation}

\begin{lemma} \label{Gamma}
Let $V$ be an algebraic finite dimensional representation over $F$
of a reductive group $G$. Then $\Gamma(V) = S_0$.
\end{lemma}
This lemma follows from fact A on page 108 in \cite{RR} for
non-archimedean $F$ and Theorem 5.2 on page 459 in \cite{Brk}.

\begin{proposition}
Let a reductive group $G$ act on an affine variety $X$. Let $x,z
\in X(F)$ be $G$-semisimple elements with different orbits. Then
there exist disjoint $G(F)$-invariant open neighborhoods $U_x$ of
$x$ and $U_z$ of $z$.
\end{proposition}
For proof of this proposition see \cite{Lun2} for archimedean $F$ and \cite{RR}, fact B on page 109 for non-archimedean $F$.

\begin{corollary} \label{EquivClassClosed}
Let a reductive group $G$ act on an affine variety $X$. Let an
element $x \in X(F)$ be $G$-semisimple. Then the set $S_x$ is
closed.
\end{corollary}
\begin{proof}
Let $y \in \overline{S_x}$. By proposition \ref{LocClosedOrbit},
$\overline{G(F)y}$ contains a closed orbit $G(F)z$. If $G(F)z = G(F)x$ then $y \in S_x$.

Otherwise, choose disjoint open $G$-invariant neighborhoods $U_z$
of $z$ and $U_x$ of $x$. Since $z \in \overline{G(F)y}$, $U_z$
intersects $G(F)y$ and hence includes $y$. Since  $y \in
\overline{S_x}$, this means that $U_z$ intersects $S_x$. Let $t
\in U_z \cap S_x$. Since $U_z$ is $G(F)$-invariant, $G(F)t \subset
U_z$. By the definition of $S_x$, $x \in \overline{G(F)t}$ and
hence $x \in \overline{U_z}$. Hence $U_z$ intersects $U_x$ -
contradiction!
\end{proof}

\subsubsection{Analytic Luna slice}

\begin{definition}
Let a reductive group $G$ act on an affine variety $X$. Let $\pi:
X(F) \to X/G(F)$ be the standard projection. An open subset $U
\subset X(F)$ is called \textbf{saturated} if there exists an open
subset $V \subset X/G(F)$ such that $U = \pi^{-1}(V)$.
\end{definition}

\noindent We will use the following corollary from Luna slice
theorem (for proof see Appendix \ref{AppLun}):

\begin{theorem} \label{LocLuna}
Let a reductive group $G$ act on a smooth affine variety $X$. Let
$x \in X(F)$ be $G$-semisimple. Then there exist\\
(i) an open $G(F)$-invariant $B$-analytic neighborhood $U$ of
$G(F)x$ in $X(F)$ with a
$G$-equivariant $B$-analytic retract $p:U \to G(F)x$ and\\
(ii) a $G_x$-equivariant $B$-analytic embedding $\psi:p^{-1}(x)
\hookrightarrow N_{Gx,x}^{X}(F)$ with open saturated image such
that $\psi(x)=0$.
\end{theorem}

\begin{definition}
In the notations of the previous theorem, denote $S:= p^{-1}(x)$
and $N:=N_{Gx,x}^{X}(F)$. We call the quintet $(U,p,\psi,S,N)$ an
\textbf{analytic Luna slice at $x$}.
\end{definition}

\begin{corollary} \label{LocLunCor}
In the notations of the previous theorem, let $y\in p^{-1}(x)$.
Denote $z:=\psi(y)$. Then\\
(i) $(G(F)_x)_z=G(F)_y$\\
(ii) $N_{G(F)y,y}^{X(F)} \cong N_{G(F)_x z, z}^{N}$ as
$G(F)_y$-spaces\\
(iii) $y$ is $G$-semisimple if and only if $z$ is
$G_x$-semisimple.
\end{corollary}

\subsection{Vector systems}$ $\\
In this subsection we introduce the term "vector system". This
term allows to formulate statements in wider generality. However,
often this generality is not necessary and therefore the reader
can skip this subsection and ignore vector systems during the
first reading.

\begin{definition}
For an analytic manifold $M$ we define the notions of
\textbf{vector system} and \textbf{B-vector system} over it.

For a smooth manifold $M$,  a vector system over $M$ is a pair
$(E,B)$ where $B$ is a smooth locally trivial fibration over $M$
and $E$ is a smooth vector bundle over $B$.

For a Nash manifold $M$, a B-vector system over $M$ is a pair
$(E,B)$ where $B$ is a Nash fibration over $M$ and $E$ is a Nash
vector bundle over $B$.

For an $l$-space $M$, a vector system over $M$ (or a B-vector
system over $M$) is an $l$-sheaf, that is locally constant sheaf,
over $M$.
\end{definition}

\begin{definition}
Let $\mathcal{V}$ be a vector system over a point $pt$. Let $M$ be
an analytic manifold. A \textbf{constant vector system with fiber
$\mathcal{V}$} is the pullback of $\mathcal{V}$ with respect to
the map $M \to pt$. We denote it by $\mathcal{V}_M$.
\end{definition}

\subsection{Preliminaries on distributions} \label{PrelDist}

\begin{definition}
Let $M$ be an analytic manifold over $F$. We define
$C_c^{\infty}(M)$ in the following way.

If $F$ is non-archimedean, $C_c^{\infty}(M)$ is the space of
locally constant compactly supported complex valued functions on
$M$. We consider no topology on it.

If $F$ is archimedean, $C_c^{\infty}(M)$ is the space of smooth
compactly supported complex valued functions on $M$. We consider
the standard topology on it.

For any analytic manifold $M$, we define the space of
distributions $\cD(M)$ by $\cD(M):=C_c^{\infty}(M)^*$. We consider
the weak topology on it.
\end{definition}

\begin{definition}
Let $M$ be a $B$-analytic manifold. We define $\Sc(M)$ in the
following way.

If $M$ is an analytic manifold over non-archimedean field,
$\Sc(M):=C_c^{\infty}(M)$.

If $M$ is a Nash manifold, $\Sc(M)$ is the space of Schwartz
functions on $M$. Schwartz functions are smooth functions that
decrease rapidly together with all their derivatives. For precise
definition see \cite{AG1}. We consider $\Sc(M)$ as a \Fre space.

For any $B$-analytic manifold $M$, we define the space of Schwartz
distributions $\Sc^*(M)$ by $\Sc^*(M):=\Sc(M)^*$.
\end{definition}

\begin{definition}
Let $M$ be an analytic manifold and let $N \subset M$ be a closed
subset. We denote
$$\cD_M(N):= \{\xi \in \cD(M)|\Supp(\xi) \subset N\}.$$

For locally closed subset $N \subset M$ we denote
$\cD_M(N):=\cD_{M\setminus (\overline{N} \setminus N)}(N)$.

Similarly we introduce the notation $\Sc^*_N(M)$ for a B-analytic
manifold $M$.
\end{definition}

\begin{definition}
Let $M$ be an analytic manifold over $F$ and $\E$ be a vector
system over $M$. We define $C_c^{\infty}(M,\E)$ in the following
way:

If $F$ is non-archimedean then $C_c^{\infty}(M,\E)$ is the space
of compactly supported sections of $\E$.

If $F$ is archimedean and $\E = (E,B)$ where $B$ is a fibration
over $M$ and $E$ is a vector bundle over $B$, then
$C_c^{\infty}(M,\E)$ is the complexification of the space of
smooth compactly supported sections of $E$ over $B$.

If $\mathcal{V}$ is a vector system over a point, we denote
$C_c^{\infty}(M,\mathcal{V}) := C_c^{\infty}(M,\mathcal{V}_M)$.
\end{definition}

We define $\cD(M,\E)$, $\cD_M(N,\E)$, $\Sc(M,\E)$, $\Sc^*(M,\E)$
and $\Sc^*_M(N,\E)$ in the natural way.

\begin{theorem}\label{Filt_nonarch}
Let an $l$-group $K$ act on an $l$-space $M$. Let $M =
\bigcup_{i=0}^l M_i$ be a $K$-invariant stratification of $M$. Let
$\chi$ be a character of $K$. Suppose that
$\Sc^*(M_i)^{K,\chi}=0$. Then $\Sc^*(M)^{K,\chi}=0$.
\end{theorem}
This theorem is a direct corollary from corollary 1.9 in
\cite{BZ}.


For the proof of the next theorem see e.g. \cite[\S B.2]{AGS1}.

\begin{theorem} \label{NashFilt}
Let a Nash group $K$ act on a Nash manifold $M$. Let $N$ be a
locally closed subset. Let $N = \bigcup_{i=0}^l N_i$ be a Nash
$K$-invariant stratification of $N$. Let $\chi$ be a character of
$K$. Suppose that for any $k \in \Z_{\geq 0}$ and $0 \leq i \leq
l$, $$\Sc^*(N_i,\Sym^k(CN_{N_i}^M))^{K,\chi}=0.$$ Then
$\Sc^*_M(N)^{K,\chi}=0.$
\end{theorem}

\begin{theorem}[Frobenius reciprocity] \label{Frob}
Let an analytic group $K$ act on an analytic manifold $M$. Let $N$
be a $K$-transitive analytic manifold. Let $\phi:M \to N$ be a
$K$-equivariant map.

Let $z \in N$ be a point and $M_z:= \phi^{-1}(z)$ be its fiber.
Let $K_z$ be the stabilizer of $z$ in $K$. Let $\Delta_K$ and
$\Delta_{K_z}$ be the modular characters of $K$ and $K_z$.

Let $\E$ be a $K$-equivariant vector system over $M$. Then\\
(i) there exists a canonical isomorphism $$Fr: \cD(M_z,\E|_{M_z}
\otimes \Delta_K|_{K_z} \cdot \Delta_{K_z}^{-1})^{K_z} \cong
\cD(M,\E)^K.$$ In particular, $Fr$ commutes with restrictions to
open sets.

(ii) For B-analytic manifolds $Fr$ maps $\Sc^*(M_z,\E|_{M_z}
\otimes \Delta_K|_{K_z} \cdot \Delta_{K_z}^{-1})^{K_z}$ to
$\Sc^*(M,\E)^K$.
\end{theorem}

For proof of (i) see \cite{Ber} 1.5 and \cite{BZ} 2.21 - 2.36 for
the case of $l$-spaces and theorem 4.2.3 in \cite{AGS1} or
\cite{Bar} for smooth manifolds. For proof of (ii) see Appendix
\ref{AppSubFrob}.

We will also use the following straightforward proposition.

\begin{proposition} \label{Product}
Let $\Omega_i \subset K_i$ be analytic groups acting on analytic
manifolds $M_i$ for $i=1 \ldots n$. Let $\E_i \to M_i$ be
$K_i$-equivariant vector systems. Suppose that
$\cD(M_i,E_i)^{\Omega_i}=\cD(M_i,E_i)^{K_i}$ for all $i$. Then
$$\cD(\prod M_i, \boxtimes E_i)^{\prod \Omega_i}=\cD(\prod M_i,
\boxtimes E_i)^{\prod K_i},$$ where $\boxtimes$ denotes the
external product.

Moreover, if $\Omega_i$, $K_i$, $M_i$ and $\E_i$ are $B$-analytic
then the same statement holds for Schwartz distributions.
\end{proposition}

For proof see e.g. \cite{AGS1}, proof of Proposition 3.1.5.

\section{Generalized Harish-Chandra descent} \label{SecDescent}

\subsection{Generalized Harish-Chandra descent}$ $\\ In this subsection
we will prove the following theorem.

\begin{theorem} \label{Gen_HC}
Let a reductive group $G$ act on a smooth affine variety $X$. Let
$\chi$ be a character of $G(F)$. Suppose that for any
$G$-semisimple $x \in X(F)$ we have
$$\cD(N_{Gx,x}^X(F))^{G(F)_x,\chi}=0.$$ Then
$$\cD(X(F))^{G(F),\chi}=0.$$
\end{theorem}

\begin{remark}
In fact, the converse is also true. We will not prove it since we
will not use it.
\end{remark}

For the proof of this theorem we will need the following lemma

\begin{lemma}
Let a reductive group $G$ act on a smooth affine variety $X$. Let
$\chi$ be a character of $G(F)$. Let $U \subset X(F)$ be an open
saturated subset. Suppose that $\cD(X(F))^{G(F),\chi}=0.$ Then
$\cD(U)^{G(F),\chi}=0.$
\end{lemma}
\begin{proof}

Consider the quotient $X/G$. It is an affine algebraic variety.
Embed it to an affine space $\mathbb{A}^n$. This defines a map
$\pi:X(F) \to F^n$. Let $V \subset X/G(F)$ be an open subset such
that $U=\pi^{-1}(V)$. There exists an open subset $V' \subset F^n$
such that $V' \cap X/G(F) = V$.

Let $\xi \in \cD(U)^{G(F),\chi}$. Suppose that $\xi$ is non-zero.
Let $x \in \Supp \xi$ and let $y:=\pi(x)$. Let $g \in
C^{\infty}_c(V')$ be such that $g(y) = 1$. Consider $\xi' \in
\cD(X(F))$ defined by $\xi'(f) := \xi(f \cdot (g \circ \pi))$.
Clearly, $x \in \Supp(\xi')$ and $\xi' \in \cD(X(F))^{G(F),\chi}$.
Contradiction.
\end{proof}

\begin{proof}[Proof of the theorem.]
Choose a $G$-semisimple $x \in X(F)$. Let ($U_x$,$p_x$,$\psi_x$,
$S_x$,$N_x$) be an analytic Luna slice at $x$.

Let $\xi' = \xi |_{U_x}$. Then $\xi' \in \cD(U_x)^{G(F),\chi}.$ By
Frobenius reciprocity it corresponds to $\xi'' \in
\cD(S_x)^{G_x(F),\chi}$.

The distribution $\xi''$ corresponds to a distribution $\xi''' \in
\cD(\psi_x(S_x))^{G_x(F),\chi}.$

However, by the previous lemma the assumption implies that
$\cD(\psi_x(S_x))^{G_x(F),\chi}=0.$ Hence $\xi'=0$.

Let $S \subset X(F)$ be the set of all $G$-semisimple points. Let
$U= \bigcup_{x\in S} U_x$. We saw that $\xi|_U=0$. On the other
hand, $U$ includes all the closed orbits, and hence by Proposition
\ref{OpenClosedAll} $U=X$.
\end{proof}

The following generalization of this theorem is proven in the same
way.

\begin{theorem} \label{Gen_HC_K}
Let a reductive group $G$ act on a smooth affine variety $X$. Let
$K \subset G(F)$ be an open subgroup and let $\chi$ be a character
of $K$. Suppose that for any $G$-semisimple $x \in X(F)$ we have
$$\cD(N_{Gx,x}^X(F))^{K_x,\chi}=0.$$ Then
$$\cD(X(F))^{K,\chi}=0.$$
\end{theorem}

Now we would like to formulate a slightly more general version of
this theorem concerning $K$-equivariant vector systems. During
first reading of this paper one can skip to the next subsection.

\begin{definition}
Let a reductive group $G$ act on a smooth affine variety $X$. Let
$K \subset G(F)$ be an open subgroup. Let $\E$ be a
$K$-equivariant vector system on $X(F)$. Let $x \in X(F)$ be
$G$-semisimple. Let $\E'$ be a $K_x$-equivariant vector system on
$N_{Gx,x}^X(F)$. We say that $\E$ and $\E'$ are
\textbf{compatible} if there exists an analytic Luna slice
$(U,p,\psi,S,N)$ such that $\E|_S = \psi^*(\E')$.
\end{definition}
Note that if $\E$ and $\E'$ are constant with the same fiber then
they are compatible.

The following theorem is proven in the same way as Theorem
\ref{Gen_HC}.
\begin{theorem} \label{Gen_HC_Sys}
Let a reductive group $G$ act on a smooth affine variety $X$. Let
$K \subset G(F)$ be an open subgroup and let $\E$ be a
$K$-equivariant vector system on $X(F)$. Suppose that for any
$G$-semisimple $x \in X(F)$ there exists a $K$-equivariant vector
system $\E'$ on $N_{Gx,x}^X(F)$, compatible with $\E$ such that
$$\cD(N_{Gx,x}^X(F),\E')^{K_x}=0.$$ Then
$$\cD(X(F),\E)^{K}=0.$$
\end{theorem}

If $\E$ and $\E'$ are B-vector systems and $K$ is open B-analytic
subgroup\footnote{In fact, any open subgroup of a B-analytic group
is B-analytic.} then the theorem holds also for Schwartz
distributions. Namely, if $\Sc^*(N_{Gx,x}^X(F),\E')^{K_x}=0$ for
any $x$ then $\Sc^*(X(F),\E)^{K}=0$, and the proof is the same.

\subsection{A stronger version}
$ $\\
In this section we give a way to validate the conditions of
theorems \ref{Gen_HC}, \ref{Gen_HC_K} and \ref{Gen_HC_Sys} by
induction.

The goal of this section is to prove the following theorem.

\begin{theorem} \label{Strong_HC}
Let a reductive group $G$ act on a smooth affine variety $X$. Let
$K \subset G(F)$ be an open subgroup and let $\chi$ be a character
of $K$. Suppose that for any $G$-semisimple $x \in X(F)$ such that
$$\cD(R(N_{Gx,x}^X))^{K_x,\chi}=0$$ we have
$$\cD(Q(N_{Gx,x}^X))^{K_x,\chi}=0.$$ Then for any for any $G$-semisimple $x \in X(F)$ we have $$\cD(N_{Gx,x}^X(F))^{K_x,\chi}=0.$$
\end{theorem}

This theorem together with Theorem \ref{Gen_HC_K} give the
following corollary.
\begin{corollary} \label{Strong_HC_Cor}
Let a reductive group $G$ act on a smooth affine variety $X$. Let
$K \subset G(F)$ be an open subgroup and let $\chi$ be a character
of $K$. Suppose that for any $G$-semisimple $x \in X(F)$ such that
$$\cD(R(N_{Gx,x}^X))^{K_x,\chi}=0$$ we have
$$\cD(Q(N_{Gx,x}^X))^{K_x,\chi}=0.$$ Then $\cD(X(F))^{K,\chi}=0.$
\end{corollary}

From now till the end of the section we fix $G$, $X$, $K$ and
$\chi$. Let us introduce several definitions and notations.

\begin{notation}
Denote \itemize  \item  $T \subset X(F)$ the set of all
$G$-semisimple points.
\item For $x,y \in T$ we say that $x>y$ if $G_x\supsetneqq G_y$.
\item $T_0 := \{x \in T | \cD(Q(N_{Gx,x}^X))^{K_x,\chi}=0 \}.$
\end{notation}
Note that if $x \in T_0$ then $\cD(N_{Gx,x}^X(F))^{K_x,\chi}=0$.

\begin{proof}[Proof of Theorem \ref{Strong_HC}]
We have to show that $T=T_0$. Assume the contrary.

Note that every chain in $T$ with respect to our ordering has a
minimum. Hence by Zorn's lemma every non-empty set in $T$ has a
minimal element. Let $x$ be a minimal element of $T - T_0$. To get
a contradiction, it is enough to show that
$\cD(R(N_{Gx,x}^X))^{K_x,\chi}=0$.

Denote $R:=R(N_{Gx,x}^X)$. By Theorem \ref{Gen_HC_K}, it is enough
to show that for any $y\in R$ we have
$$\cD(N_{G(F)_x y,y}^R)^{(K_x)_y,\chi}=0.$$

Let $(U,p,\psi,S,N)$ be an analytic Luna slice at $x$.

We can assume that $y \in \psi(S)$ since $\psi(S)$ is open,
includes 0, and we can replace $y$ by $\lambda y$ for any $\lambda
\in F^{\times}$. Let $z \in S$ be such that $\psi(z)=y$. By
corollary \ref{LocLunCor}, $(G(F)_x)_y = G(F)_z$ and $N_{G(F)_x
y,y}^R \cong N_{G z,z}^X(F)$. Hence $(K_x)_y = K_z$ and therefore
$$\cD(N_{G(F)_x y,y}^R)^{(K_x)_y,\chi} \cong \cD(N_{G
z,z}^X(F))^{K_z,\chi}.$$

However $z < x$ and hence $z \in T_0$ which means $\cD(N_{G
z,z}^X(F))^{K_z,\chi}=0$.
\end{proof}

\begin{remark}
As before, Theorem \ref{Strong_HC} and Corollary
\ref{Strong_HC_Sys} hold also for Schwartz distributions, and the
proof is the same.
\end{remark}

Again, we can formulate a more general version of Corollary
\ref{Strong_HC_Cor} concerning vector systems. During first
reading of this paper one can skip to the next subsection.

\begin{theorem} \label{Strong_HC_Sys}
Let a reductive group $G$ act on a smooth affine variety $X$. Let
$K \subset G(F)$ be an open subgroup and let $\E$ be a
$K$-equivariant vector system on $X(F)$.

Suppose that for any
$G$-semisimple $x \in X(F)$ such that\\
(*) any $K_x \times F^{\times}$-equivariant vector system $\E'$ on
$R(N_{Gx,x}^X)$ compatible with $\E$ satisfies
$\cD(R(N_{Gx,x}^X),\E')^{K_x}=0$ (where the action of $F^{\times}$
is the homothety action),

we have\\
(**)  there exists a $K_x \times F^{\times}$-equivariant vector
system $\E'$ on $Q(N_{Gx,x}^X)$ compatible with $\E$ such that
$$\cD(Q(N_{Gx,x}^X),\E')^{K_x}=0.$$

Then $\cD(X(F),\E)^{K}=0.$
\end{theorem}

The proof is the same as the proof of Theorem \ref{Strong_HC}
using the following lemma that follows from the definitions.

\begin{lemma}
Let a reductive group $G$ act on a smooth affine variety $X$. Let
$K \subset G(F)$ be an open subgroup and let $\E$ be a
$K$-equivariant vector system on $X(F)$. Let $x \in X(F)$ be
$G$-semisimple. Let $(U, p, \psi, S, N)$ be an analytic Luna slice
at $x$.

Let $\E'$ be a $K_x$-equivariant vector system on $N$ compatible
with $\E$. Let $y \in S$ be $G$-semisimple. Let $z := \psi(y)$.
Let $\E''$ be a $(K_x)_z$-equivariant vector system on
$N_{G_xz,z}^N$ compatible with $\E'$. Consider the isomorphism
$N_{G_xz,z}^N(F) \cong N_{Gy,y}^X(F)$ and let $\E'''$ be the
corresponding $K_y$-equivariant vector system on $N_{Gy,y}^X(F)$.

Then $\E'''$ is compatible with $\E$.
\end{lemma}

Again, if $\E$ and $\E'$ are B-vector systems then the theorem
holds also for Schwartz distributions.

\section{Distributions versus Schwartz distributions}
\label{DistVerSch}

\setcounter{lemma}{0}

The tools developed in the previous section enabled us to prove
the following version of localization principle.

\begin{theorem}[Localization principle] \label{LocPrin}
Let a reductive group $G$ act on a smooth algebraic variety $X$.
Let $Y$ be an algebraic variety and $\phi:X \to Y$ be an affine
algebraic $G$-invariant map. Let $\chi$ be a character of $G(F)$.
Suppose that for any $y \in Y(F)$ we have
$\cD_{X(F)}(\phi(F)^{-1}(y))^{G(F),\chi}=0$. Then
$\cD(X(F))^{G(F),\chi}=0$.
\end{theorem}

For proof see Appendix
\ref{SecRedLocPrin}.

In this section we use this theorem to show that if there are no
$G(F)$-equivariant Schwartz distributions on $X(F)$ then there are
no $G(F)$-equivariant distributions on $X(F)$.

\begin{theorem}   \label{NoSNoDist}
Let a reductive group $G$ act on a smooth affine variety $X$. Let
$V$ be a finite dimensional continuous representation of $G(F)$
over $\R$. Suppose that $\Sc^*(X(F),V)^{G(F)}=0$. Then
$\cD(X(F),V)^{G(F)}=0$.
\end{theorem}

For the proof we will need the following definition and theorem.

\begin{definition} $ $

(i) Let a topological group $K$ act on a topological
space $M$. We call a closed $K$-invariant subset $C \subset M$ \textbf{compact
modulo $K$} if there exists a compact subset $C' \subset M$ such
that $C \subset KC'.$

(ii) Let a Nash group $K$ act on a Nash
manifold $M$. We call a closed $K$-invariant subset $C \subset M$ \textbf{Nashly compact
modulo $K$} if there exist a compact subset $C' \subset M$ and semi-algebraic closed subset $Z \subset M$ such
that $C \subset Z \subset KC'.$
\end{definition}

\begin{remark}
Let a reductive group $G$ act on a smooth affine variety $X$. Let $K:=G(F)$ and $M:=X(F)$. Then it is easy to see that the notions of
compact modulo $K$ and Nashly compact modulo $K$ coincide.
\end{remark}

\begin{theorem} \label{CompSchwartz}
Let a Nash group $K$ act on a Nash manifold $M$. Let $E$ be a
$K$-equivariant Nash bundle over $M$. Let $\xi \in \cD(M,E)^K$
such that $\Supp(\xi)$ is Nashly compact modulo $K$. Then $\xi \in
\Sc^*(M,E)^K$.
\end{theorem}

The formulation and the idea of the proof of this theorem are due
to J. Bernstein. For the proof see Appendix \ref{KInvAreSchwartz}.

\begin{proof}[Proof of Theorem  \ref{NoSNoDist}]
Fix any $y \in X/G(F)$ and denote $M:=\pi_X^{-1}(y)(F)$.

By localization principle (Theorem \ref{LocPrin} and Remark
\ref{RemLocVectSys}), it is enough to prove that
$$\Sc^*_{X(F)}(M,V)^{G(F)}=\cD_{X(F)}(M,V)^{G(F)}.$$ Choose $\xi
\in \cD_{X(F)}(M,V)^{G(F)}$. $M$ has a unique stable closed
$G$-orbit and hence a finite number of closed $G(F)$-orbits. By
Theorem \ref{CompSchwartz}, it is enough to show that $M$ is Nashly
compact modulo $G(F)$. Clearly $M$ is semi-algebraic. Choose representatives $x_i$ of the closed
$G(F)$ orbits in $M$. Choose compact neighborhoods $C_i$ of $x_i$.
Let $C' := \bigcup C_i$. By corollary \ref{OpenClosedAll}, $G(F)C'
\supset M$.
\end{proof}

\section{Applications of Fourier transform and Weil
representation} \label{SecFour}

Let $G$ be a reductive group. Let $V$ be a finite dimensional
algebraic representation of $G$ over $F$. Let $\chi$ be a
character of $G(F)$. In this section we provide some tools to
verify $\Sc^*(Q(V))^{G(F),\chi}=0$ if we know that
$\Sc^*(R(V))^{G(F),\chi}=0$.

\subsection{Preliminaries}
$ $\\
From now till the end of the paper we fix an additive character
$\kappa$ of $F$. If $F$ is archimedean we fix $\kappa$ to be
defined by $\kappa(x):=e^{2\pi i \re(x)}$.

\begin{notation}
Let $V$ be a vector space over $F$. Let $B$ be a non-degenerate
bilinear form on $V$. We denote by $\Fou_B: \Sc^*(V) \to \Sc^*(V)$
the Fourier transform given by $B$ with respect to the
self-adjoint Haar measure on $V$. For any B-analytic manifold $M$
over $F$ we also denote by $\Fou_B:\Sc^*(M \times V) \to \Sc^*(M
\times V)$ the partial Fourier transform.
\end{notation}

\begin{notation}
Let $V$ be a vector space over $F$. Consider the homothety action
of $F^{\times}$ on $V$ by $\rho(\lambda)v:= \lambda^{-1}v$. It
gives rise to an action $\rho$ of $F^{\times}$ on $\Sc^*(V)$.

Also, for any $\lambda \in F^{\times}$ denote
$|\lambda|:=\frac{dx}{\rho(\lambda)dx}$, where $dx$ denotes the
Haar measure on $F$. Note that for $F =\R$, $|\lambda|$ is equal
to the classical absolute value but for $F =\C$, $|\lambda| = (Re
\lambda)^2+ (Im \lambda)^2$.
\end{notation}

\begin{notation}
Let $V$ be a vector space over $F$. Let $B$ be a non-degenerate
symmetric bilinear form on $V$. We denote by $\gamma(B)$ the Weil
constant.  For its definition see e.g. \cite{Gel}, section 2.3 for
non-archimedean $F$ and \cite{RS1}, section 1 for archimedean $F$.

For any $t\in F^{\times}$ denote $\delta_B(t)=
\gamma(B)/\gamma(tB)$.
\end{notation}

Note that $\gamma_B(t)$ is an eights root of unity and if $dimV$
is odd and $F \neq \C$ then $\delta_B$ is \textbf{not} a
multiplicative character.

\begin{notation}
Let $V$ be a vector space over $F$. Let $B$ be a non-degenerate
symmetric bilinear form on $V$. We denote $$Z(B):=\{x \in
V|B(x,x)=0 \}.$$
\end{notation}

\begin{theorem} [non-archimedean homogeneity] \label{NonArchHom}
Suppose that $F$ is \textbf{non-archimedean}. Let $V$ be a vector
space over $F$. Let $B$ be a non-degenerate symmetric bilinear
form on $V$. Let $M$ be a B-analytic manifold over $F$. Let $\xi
\in \Sc^*_{V\times M}(Z(B)\times M)$ be such that $\Fou_B(\xi) \in
\Sc^*_{V\times M}(Z(B)\times M)$. Then for any $t \in F^{\times}$,
we have $\rho(t)\xi=\delta_B(t) |t|^{dimV/2} \xi$ and $\xi=
\gamma(B)^{-1} \Fou_B \xi$. In particular, if $dimV$ is odd then
$\xi = 0$.
\end{theorem}
For proof see \cite{RS2}, section 8.1.

For the archimedean version of this theorem we will need the
following definition.
\begin{definition}
Let $V$ be a finite dimensional vector space over $F$. Let $B$ be
a non-degenerate symmetric bilinear form on $V$. Let $M$ be a
B-analytic manifold over $F$. We say that a distribution $\xi \in
\Sc^*(V \times M)$ is
\textbf{adapted to $B$} if either \\
(i) for any $t \in F^{\times}$ we have $\rho(t)\xi=
\delta(t)|t|^{dimV/2} \xi$ and $ \xi$ is proportional to $\Fou_B \xi$ or\\
(ii) $F$ is archimedean and for any $t \in F^{\times}$ we have
$\rho(t)\xi= \delta(t)t|t|^{dimV/2} \xi$.
\end{definition}

Note that if $dimV$ is odd and $F \neq \C$ then every $B$-adapted
distribution is zero.

\begin{theorem} [archimedean homogeneity] \label{ArchHom}
Let $V$ be a vector space over $F$. Let $B$ be a non-degenerate
symmetric bilinear form on $V$. Let $M$ be a Nash manifold. Let $L
\subset \Sc^*_{V\times M}(Z(B)\times M)$ be a non-zero subspace
such that $\forall \xi \in L $ we have $\Fou_B(\xi) \in L$ and $B
\xi \in L$ (here $B$ is interpreted as a quadratic form).

Then there exists a non-zero distribution $\xi \in L$ which is
adapted to $B$.
\end{theorem}
For archimedean $F$ we prove this theorem in Appendix
\ref{AppRealHom}. For the non-archimedean $F$ it follows from
Theorem \ref{NonArchHom}.

We will also use the following trivial observation.

\begin{lemma}
Let $V$ be a finite dimensional vector space over $F$. Let a
B-analytic group $K$ act linearly on $V$. Let $B$ be a
$K$-invariant non-degenerate symmetric bilinear form on $V$. Let
$M$ be a B-analytic $K$-manifold over $F$. Let $\xi \in \Sc^*(V
\times M)$ be a $K$-invariant distribution. Then $\Fou_B(\xi)$ is
also $K$-invariant.
\end{lemma}

\subsection{Applications}
$ $\\
The following two theorems easily follow form the results of the
previous subsection.

\begin{theorem} \label{Non_arch_Homog}
Suppose that $F$ is non-archimedean. Let $G$ be a reductive group.
Let $V$ be a finite dimensional algebraic representation of $G$
over $F$. Let $\chi$ be character of $G(F)$. Suppose that
$\Sc^*(R(V))^{G(F),\chi}=0$. Let $V = V_1 \oplus V_2$ be a
$G$-invariant decomposition of $V$. Let $B$ be a $G$-invariant
symmetric non-degenerate bilinear form on $V_1$. Consider the
action $\rho$ of $F^{\times}$ on $V$ by homothety on $V_1$.

Then any $\xi \in \Sc^*(Q(V))^{G(F),\chi}$ satisfies
$\rho(t)\xi=\delta_B(t) |t|^{dimV_1/2} \xi$ and $\xi= \gamma(B)
\Fou_B \xi$. In particular, if $dimV_1$ is odd then $\xi = 0$.
\end{theorem}

\begin{theorem} \label{Homog}
Let $G$ be a reductive group. Let $V$ be a finite dimensional
algebraic representation of $G$ over $F$. Let $\chi$ be character
of $G(F)$. Suppose that $\Sc^*(R(V))^{G(F),\chi}=0$. Let $Q(V) = W
\oplus (\bigoplus_{i=1}^kV_i)$ be a $G$-invariant decomposition of
$Q(V)$. Let $B_i$ be $G$-invariant symmetric non-degenerate
bilinear forms on $V_i$. Suppose that any $\xi \in
\Sc^*_{Q(V)}(\Gamma(V))^{G(F),\chi}$ which is adapted to each
$B_i$ is zero.

Then $\Sc^*(Q(V))^{G(F),\chi}=0.$
\end{theorem}

\begin{remark}
One can easily generalize theorems \ref{Homog} and
\ref{Non_arch_Homog} to the case of constant vector systems.
\end{remark}

\section{Tame actions} \label{SecTame}
\setcounter{lemma}{0}

In this section we consider problems of the following type. A
reductive group $G$ acts on a smooth affine variety $X$, and
$\tau$ is an automorphism of $X$ which normalizes the action of
$G$. We want to check whether any $G(F)$-invariant Schwartz
distribution on $X(F)$ is also $\tau$-invariant.

\begin{definition}
Let $\pi$ be an action of a reductive group $G$ on a smooth affine variety $X$.
We say that an algebraic automorphism $\tau$ of $X$ is \textbf{$G$-admissible} if \\
(i) $\pi(G(F))$ is of index $\leq 2$ in the group of automorphisms
of $X$
generated by $\pi(G(F))$ and $\tau$.\\
(ii) For any closed $G(F)$ orbit $O \subset X(F)$, we have
$\tau(O)=O$.
\end{definition}

\begin{proposition} \label{AdmisDescends}
Let $\pi$ be an action of a reductive group $G$ on a smooth affine
variety $X$. Let $\tau$ be a $G$-admissible automorphism of $X$.
Let $K:=\pi(G(F))$ and let $\widetilde{K}$ be the group generated
by $\pi(G(F))$ and $\tau$. Let $x \in X(F)$ be a point with closed
$G(F)$ orbit. Let $\tau' \in \widetilde{K}_x - K_x$. Then $d\tau'
|_{N_{Gx,x}^X}$ is $G_x$-admissible.
\end{proposition}

\begin{proof}
Let $\widetilde{G}$ denote the group generated by $G$ and
$\tau$.\\
(i) is obvious.\\
(ii) Let $y \in N_{Gx,x}^X(F)$ be an element with closed $G_x$
orbit. Let $y' = d\tau'(y)$. We have to show that there exists $g
\in G_x(F)$ such that $gy = gy'$. Let $(U,p,\psi,S,N)$ be analytic
Luna slice at $x$ with respect to the action of $\widetilde{G}$.
We can assume that there exists $z \in S$ such that $y=\psi(z)$.
Let $z' = \tau'(z)$. By corollary \ref{LocLunCor}, $z$ is
$G$-semisimple. Since $\tau$ is admissible, this implies that
there exists $g \in G(F)$ such that $gz = z'$. Clearly, $g \in
G_x(F)$ and $gy = y'$.
\end{proof}

\begin{definition}
We call an action of a reductive group $G$ on a smooth affine
variety $X$ \textbf{tame} if for any $G$-admissible $\tau : X \to
X$, we have $\Sc^*(X(F))^{G(F)} \subset \Sc^*(X(F))^{\tau}.$
\end{definition}

\begin{definition}
We call an algebraic representation  of a reductive group $G$ on a
finite dimensional linear space $V$ over $F$ \textbf{linearly
tame} if for any $G$-admissible linear map $\tau : V \to V$, we
have $\Sc^*(V(F))^{G(F)} \subset \Sc^*(V(F))^{\tau}.$

We call a representation \textbf{weakly linearly tame} if for any
$G$-admissible linear map $\tau : V \to V$, such that
$\Sc^*(R(V))^{G(F)} \subset \Sc^*(R(V))^{\tau}$ we have
$\Sc^*(Q(V))^{G(F)} \subset \Sc^*(Q(V))^{\tau}.$
\end{definition}

\begin{theorem} \label{Invol_HC}
Let a reductive group $G$ act on a smooth affine variety $X$.
Suppose that for any $G$-semisimple $x \in X(F)$, the action of
$G_x$ on $N_{Gx,x}^X$ is weakly linearly tame. Then the action of
$G$ on $X$ is tame.
\end{theorem}

The proof is rather straightforward except of one minor
complication: the group of automorphisms of $X(F)$ generated by
the action of $G(F)$ is not necessarily a group of $F$ points of
any algebraic group.

\begin{proof}
Let $\tau:X \to X$ be an admissible automorphism.

Let $\widetilde{G} \subset Aut(X)$ be the algebraic group
generated by the actions of $G$ and $\tau$. Let $K \subset
Aut(X(F))$ be the B-analytic group generated by the action of
$G(F)$. Let $\widetilde{K} \subset Aut(X(F))$ be the B-analytic
group generated by the actions of $G$ and $\tau$. Note that
$\widetilde{K} \subset \widetilde{G}(F)$ is an open subgroup of
finite index. Note that for any $x \in X(F)$, $x$ is
$\widetilde{G}$-semisimple if and only if it is $G$-semisimple. If
$K = \widetilde{K}$ we are done, so we will assume $K \neq
\widetilde{K}$. Let $\chi$ be the character of $\widetilde{K}$
defined by $\chi(K)=\{1\}$, $\chi(\widetilde{K} - K) = \{-1\}$.

It is enough to prove that $\Sc^*(X)^{\widetilde{K},\chi}=0$. By
generalized Harish-Chandra descent (corollary \ref{Strong_HC_Cor})
it is enough to prove that for any $G$-semisimple $x \in X$ such
that $\Sc^*(R(N_{Gx,x}^X))^{\widetilde{K}_x,\chi}=0$ we have
$\Sc^*(Q(N_{Gx,x}^X))^{\widetilde{K}_x,\chi}=0$. Choose any
automorphism $\tau'\in \widetilde{K}_x - K_x$. Note that $\tau'$
and $K_x$ generate $\widetilde{K}_x$. Denote $$\eta :=
d\tau'|_{N_{Gx,x}^X(F)}.$$ By Proposition \ref{AdmisDescends},
$\eta$ is $G$- admissible. Note that
$$\Sc^*(R(N_{Gx,x}^X))^{K_x}= \Sc^*(R(N_{Gx,x}^X))^{G(F)_x} \text{ and }
\Sc^*(Q(N_{Gx,x}^X))^{K_x}= \Sc^*(Q(N_{Gx,x}^X))^{G(F)_x}.$$
Hence we have $$\Sc^*(R(N_{Gx,x}^X))^{G(F)_x} \subset
\Sc^*(R(N_{Gx,x}^X))^{\eta}.$$ Since the action of $G_x$ is weakly
linearly tame, this implies that $$\Sc^*(Q(N_{Gx,x}^X))^{G(F)_x}
\subset \Sc^*(Q(N_{Gx,x}^X))^{\eta}$$ and therefore
$\Sc^*(Q(N_{Gx,x}^X))^{\widetilde{K}_x,\chi}=0$.
\end{proof}

\begin{definition}
We call an algebraic representation  of a reductive group $G$ on a
finite dimensional linear space $V$ over $F$ \textbf{special} if
for any $\xi \in \Sc^*_{Q(V)}(\Gamma(V))^{G(F)}$ such that for any
$G$-invariant decomposition $Q(V) = W_1 \oplus W_2$ and any two
$G$-invariant symmetric non-degenerate  bilinear forms $B_i$ on
$W_i$ the Fourier transforms $\Fou_{B_i}(\xi)$ are also supported
in $\Gamma(V)$, we have $\xi = 0$.
\end{definition}

\begin{proposition} \label{SpecWeakTameAct}
Every special algebraic representation $V$ of a reductive group
$G$ is weakly linearly tame.
\end{proposition}

This proposition follows immediately from the following lemma.

\begin{lemma}
Let $V$ be an algebraic representation of a reductive group $G$.
Let $\tau$ be an admissible linear automorphism of $V$. Let $V =
W_1 \oplus W_2$ be a $G$-invariant decomposition of $V$ and $B_i$
be $G$-invariant symmetric non-degenerate  bilinear forms on
$W_i$. Then $W_i$ and $B_i$ are also $\tau$-invariant.
\end{lemma}

This lemma follows in turn from the following one.

\begin{lemma}
Let $V$ be an algebraic representation of a reductive group $G$.
Let $\tau$ be an admissible automorphism of $V$. Then $\cO(V)^{G}
\subset \cO(V)^{\tau}$.
\end{lemma}
\begin{proof}
Consider the projection $\pi:V \to V/G$. We have to show that
$\tau$ acts trivially on $V/G$. Let $x \in \pi(V(F))$. Let
$X:=\pi^{-1}(x)$. By Proposition \ref{LocClosedOrbit} $G(F)$ has a
closed orbit in $X(F)$. The automorphism $\tau$ preserves this
orbit and hence preserves $x$. So $\tau$ acts trivially on
$\pi(V(F))$, which is Zariski dense in $V/G$. Hence $\tau$ acts
trivially on $V/G$.
\end{proof}

Now we introduce a criterion that allows to prove that a
representation is special. It follows immediately from Theorem
\ref{ArchHom}.

\begin{lemma} \label{SpecActCrit}
Let $V$ be an algebraic representation of a reductive group $G$.
Let $Q(V) = \bigoplus W_i$ be a $G$-invariant decomposition. Let
$B_i$ be symmetric non-degenerate $G$-invariant bilinear forms on
$W_i$. Suppose that any $\xi \in \Sc^*_{Q(V)}(\Gamma(V))^{G(F)}$
which is adapted to all $B_i$ is zero. Then $V$ is special.
\end{lemma}

\section{Symmetric pairs} \label{SecSymPairs}
In this section we apply our tools to symmetric pairs. We
introduce several properties of symmetric pairs and discuss their
interrelations. In Appendix \ref{Diag} we present a diagram that
illustrates the most important ones.

\subsection{Preliminaries and notations}

\begin{definition}
A \textbf{symmetric pair} is a triple $(G,H,\theta)$ where $H
\subset G$ are reductive groups, and $\theta$ is an involution of
$G$ such that $H = G^{\theta}$. We call a symmetric pair
\textbf{connected} if $G/H$ is connected.

For a symmetric pair $(G,H,\theta)$ we define an antiinvolution
$\sigma :G \to G$ by $\sigma(g):=\theta(g^{-1})$, denote $\g:=Lie
G$, $\h := LieH$. Let $\theta$ and $\sigma$ act on $\g$ by their
differentials and denote $\gd:=\{a \in \g | \sigma(a)=a\}=\{a \in
\g | \theta(a)=-a\}$. Note that $H$ acts on $\gd$ by the adjoint
action. Denote also $G^{\sigma}:=\{g \in G| \sigma(g)=g\}$ and
define a \textbf{symmetrization map} $s:G \to G^{\sigma}$ by
$s(g):=g \sigma(g)$.
\end{definition}

\begin{definition}
Let $(G_1,H_1,\theta_1)$ and $(G_2,H_2,\theta_2)$ be symmetric
pairs. We define their \textbf{product} to be the symmetric pair
$(G_1 \times G_2,H_1 \times H_2,\theta_1 \times \theta_2)$.
\end{definition}

\begin{theorem} \label{RegFun}
For any connected symmetric pair $(G,H,\theta)$ we have
$\mathcal{O}(G)^{H \times H} \subset \mathcal{O}(G)^{\sigma}$.
\end{theorem}

\begin{proof}
Consider the multiplication map $H \times G^{\sigma} \to G$. It is
\et at $1 \times 1$ and hence its image $HG^{\sigma}$ contains an
open neighborhood of 1 in $G$. Hence the image of $HG^{\sigma}$ in
$G/H$ is dense. Thus $HG^{\sigma}H$ is dense in $G$. Clearly
$\mathcal{O}(H G^{\sigma} H)^{H \times H} \subset \mathcal{O}(H
G^{\sigma} H)^{\sigma}$ and hence $\mathcal{O}(G)^{H \times H}
\subset \mathcal{O}(G)^{\sigma}$.
\end{proof}
\begin{corollary} \label{ClosedOrbits}
For any connected symmetric pair $(G,H,\theta)$ and any closed $H
\times H$ orbit $\Delta \subset G$, we have
$\sigma(\Delta)=\Delta$.
\end{corollary}
\begin{proof}
Denote $\Upsilon:=H \times H$. Consider the action of the
2-element group $(1,\tau)$ on $\Upsilon$ given by $\tau(h_1,h_2):=
(\theta(h_2), \theta(h_1))$. This defines the semi-direct product
$\widetilde{\Upsilon}:= (1,\tau) \ltimes \Upsilon $. Extend the
two-sided action of $\Upsilon$ to $\widetilde{\Upsilon}$ by the
antiinvolution $\sigma$. Note that the previous theorem implies
that $G/\Upsilon = G/\widetilde{\Upsilon}$. Let $\Delta$ be a
closed $\Upsilon$-orbit. Let $\widetilde{\Delta}:=\Delta \cup
\sigma(\Delta)$. Let $a := \pi_G(\widetilde{\Delta}) \subset
G/\widetilde{\Upsilon}$. Clearly, $a$ consists of one point. On
the other hand, $G/\widetilde{\Upsilon} = G/\Upsilon$ and hence
$\pi_G^{-1}(a)$ contains a unique closed $G$-orbit. Therefore
$\Delta = \widetilde{\Delta} = \sigma(\Delta)$.
\end{proof}

\begin{corollary}
Let $(G,H,\theta)$ be a connected symmetric pair. Let $g \in G(F)$
be $H\times H$-semisimple. Suppose that $H^1(F,(H \times H)_g)$ is
trivial. Then $\sigma(g) \in H(F)gH(F)$.
\end{corollary}
For example, if $(H \times H)_g$ is a product of general linear
groups over some field extensions then $H^1(F,(H \times H)_g)$ is
trivial.

\begin{definition}
A symmetric pair $(G,H,\theta)$ is called \textbf{good}
if for any closed $H(F) \times H(F)$ orbit $O \subset G(F)$, we
have $\sigma(O)=O$.
\end{definition}

\begin{corollary} \label{ComplexGood}
Any connected symmetric pair over $\C$ is good.
\end{corollary}

\begin{definition}
A symmetric pair $(G,H,\theta)$ is called a \textbf{GK pair} if
$$\Sc^*(G(F))^{H(F) \times H(F)} \subset \Sc^*(G(F))^{\sigma}.$$
\end{definition}

We will see later in section \ref{Gel} that GK pairs satisfy a
Gelfand pair property that we call GP2 (see Definition \ref{GPs}
and Theorem \ref{DistCrit}). Clearly every GK pair is good and we
conjecture that the converse is also true. We will discuss it in
more details in subsection \ref{conj}.

\begin{lemma} \label{BilForm}
Let $(G,H,\theta)$ be a symmetric pair. Then there exists a
$G$-invariant $\theta$-invariant non-degenerate symmetric bilinear
form $B$ on $\g$. In particular, $B|_{\h}$ and $B|_{\g^{\sigma}}$
are also non-degenerate and $\h$ is orthogonal to $\g^{\sigma}$.
\end{lemma}
\begin{proof}$ $\\
\indent
Step 1. Proof for semisimple $\g$.\\
Let $B$ be the Killing form on $\g$. Since it is non-degenerate,
it is enough to show that $\h$ is orthogonal to $\g^{\sigma}$. Let
$A \in \h$ and $B \in \g^{\sigma}$. We have to show
$tr(Ad(A)Ad(B))=0$. This follows from the fact that
$Ad(A)Ad(B)(\h) \subset \g^{\sigma}$ and $Ad(A)Ad(B)(\g^{\sigma})
\subset \h$.

Step 2. Proof in the general case.\\
Let $\g = \g' \oplus \z$ such that $\g'$ is semisimple and $\z$ is
the center. It is easy to see that this decomposition is $\theta$
invariant. Now the proposition easily follows from the previous
case.
\end{proof}
\begin{lemma} \label{ExpMap}
Let $(G,H,\theta)$ be a symmetric pair. Then there exists an
$Ad(G(F))$-equivariant and $\sigma$-equivariant map
$\mathcal{U}(G) \to \mathcal{N}(\g)$ where $\mathcal{U}(G)$ is the
set of unipotent elements in $G(F)$ and $\mathcal{N}(\g)$ is the
set of nilpotent elements in $\g(F)$.
\end{lemma}
\begin{proof}
It follows from the existence of analytic Luna slice at point $1
\in G(F)$ with respect to the action of $\widetilde{G}$ where
$\widetilde{G}$ is the group generated by $\sigma$ and the adjoint
action of $G$ on itself.
\end{proof}

\begin{lemma} \label{SL2triple}
Let $(G,H,\theta)$ be a symmetric pair. Let $x \in \gd$ be a
nilpotent element. Then there exists a group homomorphism
$\phi:SL_2 \to G$ such that $$d\phi(\begin{pmatrix}
  0 & 1 \\
  0 & 0\end{pmatrix}) = x, \quad
d\phi(\begin{pmatrix}
  0 & 0 \\
  1 & 0\end{pmatrix}) \in \gd \text{ and }\quad
\phi(\begin{pmatrix}
  t & 0 \\
  0 & t^{-1}\end{pmatrix}) \in H.$$
In particular $0 \in \overline{Ad(H)(x)}$.
\end{lemma}
This lemma was essentially proven for $F=\C$ in \cite{KR}. The
same proof works for any $F$ and we repeat it here for the
convenience of the reader.

\begin{proof}

By Jacobson-Morozov theorem (see \cite{Jac}, Chapter III, Theorems
17 and 10) we can complete $x$ to an $sl_2$-triple $(x_-,s,x)$.
Let $s':= \frac{s+\theta(s)}{2}$. It satisfies $[s',x]=2x$ and
lies in the ideal $[x,\g]$ and hence by Morozov lemma (see
\cite{Jac}, Chapter III, Lemma 7), $x$ and $s'$ can be completed
to an $sl_2$ triple $(x_-,s',x)$. Let $x'_-:=\frac{x_- -
\theta(x_-)}{2}$. Note that $(x_-',s',x)$ is also an
$sl_2$-triple. Exponentiating this $sl_2$-triple to a map $SL_2
\to G$ we get the required homomorphism.
\end{proof}

\begin{notation}
In the notations of the previous lemma we denote $$D_t(x):=
\phi(\begin{pmatrix}
  t & 0 \\
  0 & t^{-1}\end{pmatrix}) \text{ and } d(x):=d\phi(\begin{pmatrix}
  1 & 0 \\
  0 & -1\end{pmatrix}).$$
Those elements depend on the choice of $\phi$. However, whenever
we will use this notation nothing will depend on their choice.
\end{notation}

\subsection{Descendants of symmetric pairs}

\begin{proposition} \label{PropDescend}
Let  $(G,H,\theta)$ be a symmetric pair. Let $g \in G(F)$ be $H
\times H$-semisimple. Let $x=s(g)$. Then \\
(i) $x$ is semisimple.\\
(ii) Consider the adjoint action of $G$ on itself and the
two-sided action of $H \times H$ on $G$. Then $H_x \cong (H \times
H)_g$ and $(\g_x)^{\sigma} \cong N_{H g H,g}^G$ as $H_x$ spaces.
\end{proposition}
\begin{proof}
$ $\\ \indent (i) Let $x = x_sx_u$ be the Jordan decomposition of
$x$. The uniqueness of Jordan decomposition implies that both
$x_u$ and $x_s$ belong to $G^{\sigma}$. To show that $x_u=1$ it is
enough to show that $\overline{Ad(H)(x)} \ni x_s$. We will do that
in several steps.
$ $\\
\indent Step 1. Proof for the case when $x_s = 1$.\\
It follows immediately from the two previous lemmas (\ref{ExpMap}
and \ref{SL2triple}).

Step 2. Proof for the case when
$x_s \in Z(G)$.\\
This case follows from Step 1 since conjugation acts trivially on
$Z(G)$.

Step 3. Proof in the general case.\\ The statement follows from
Step 2 for the group $G_{x_s}$.

(ii) The symmetrization gives rise to an isomorphism $(H \times
H)_g \cong H_x$. Let us now prove $(\g_x)^{\sigma} \cong N_{H g
H,g}^G$. First of all, $N_{H g H,g}^G \cong \g / (\h+ Ad(g)\h).$
Let $\theta'$ be the involution of $G$ defined by $\theta'(y)= x
\theta(y) x^{-1}$. Note that $Ad(g)\h= \g^{\theta'}$. Fix a
non-degenerate $G$-invariant symmetric bilinear form $B$ on $\g$
as in Lemma \ref{BilForm}. Note that $B$ is also $\theta'$
invariant and hence $$(Ad(g)\h)^{\bot} = \{a \in \g|
\theta'(a)=-a\}.$$ Now $$N_{H g H,g}^G \cong (\h+ Ad(g)\h)^{\bot}
= \h^{\bot} \cap Ad(g)\h^{\bot} = \{a \in \g|
\theta(a)=\theta'(a)=-a\}= (\g_x)^{\sigma}.$$
\end{proof}

It is easy to see that the isomorphism $N_{H g H,g}^G \cong
(\g_x)^{\sigma}$ does not depend on the choice of $B$.

\begin{definition} \label{descendant}
In the notations of the previous proposition we will say that the
pair $(G_x,H_x.\theta|_{G_x})$ is a \textbf{descendant} of
$(G,H,\theta)$.
\end{definition}

\subsection{Tame symmetric pairs} \label{SecTamePairs}

\begin{definition}
We call a symmetric pair $(G,H,\theta)$\\
(i) \textbf{tame} if the action of $H\times H$ on $G$ is tame.\\
(ii) \textbf{linearly tame} if the action of $H$ on $\g^{\sigma}$
is linearly tame.\\
(iii) \textbf{weakly linearly tame} if the action of $H$ on
$\g^{\sigma}$ is weakly linearly tame.
\end{definition}

\begin{remark}
Evidently, any good tame symmetric pair is a GK pair.
\end{remark}

The following theorem is a direct corollary of Theorem
\ref{Invol_HC}.

\begin{theorem} \label{LinDes}
Let $(G,H,\theta)$ be a symmetric pair. Suppose that all its
descendants (including itself)  are weakly linearly tame. Then
$(G,H,\theta)$ is tame and linearly tame.
\end{theorem}

\begin{definition}
We call a symmetric pair $(G,H,\theta)$ \textbf{special} if
$g^{\sigma}$ is a special representation of $H$.
\end{definition}

\begin{proposition} \label{SpecWeakReg}
Any special symmetric pair is weakly linearly tame.
\end{proposition}
This proposition follows immediately from Proposition
\ref{SpecWeakTameAct}

\begin{proposition}
A product of special symmetric pairs is special.
\end{proposition}
The proof of this proposition is straightforward using Lemma
\ref{BilForm}.\\

Now we would like to give a criterion of speciality for symmetric
pairs.

\begin{proposition}[Speciality criterion] \label{SpecCrit}
Let $(G,H,\theta)$ be a symmetric pair. Suppose that for any
nilpotent $x \in \gd$ either\\
(i) $\tr(ad(d(x))|_{\h_x}) < dim \g^{\sigma} $
or\\
(ii) $F$ is non-archimedean and  $\tr(ad(d(x))|_{\h_x}) \neq dim
\g^{\sigma}$.

Then the pair $(G,H,\theta)$ is special.
\end{proposition}

For the proof we will need the following lemmas.

\begin{lemma} \label{EquivDef}
Let $(G,H,\theta)$ be a symmetric pair. Then $\Gamma(\gd)$ is the
set of all nilpotent elements in $Q(\gd)$.
\end{lemma}
This lemma is a direct corollary from Lemma \ref{SL2triple}.

\begin{lemma} \label{EigenInt}
Let $(G,H,\theta)$ be a symmetric pair. Let $x \in \gd$ be a
nilpotent element. Then all the eigenvalues of
$ad(d(x))|_{\gd/[x,\h]}$ are non-positive integers.
\end{lemma}

This lemma follows from the existence of a natural onto map $
\g/[x,\g] \twoheadrightarrow \gd/[x,\h]$ using the following
straightforward lemma.

\begin{lemma}
Let $V$ be a representation of an $sl_2$ triple $(e,h,f)$. Then
all the eigenvalues of $h|_{V/e(V)}$ are non-positive integers.
\end{lemma}

Now we are ready to prove the speciality criterion.
\begin{proof}[Proof of Proposition \ref{SpecCrit}]
We will give a proof in the case that $F$ is archimedean. The case
of non-archimedean $F$ is done in the same way but with less
complications.

Let $\chi$ be a character of $F^{\times}$ given by either
$\chi(\lambda)=u(\lambda)|\lambda|^{dim\gd/2}$ or
$\chi(\lambda)=u(\lambda)|\lambda|^{dim\gd/2+1}$, where $u$ is
some unitary character. By Lemma \ref{SpecActCrit} it is enough to
prove
$$\Sc^*_{Q(\gd)}(\Gamma(\gd))^{H(F) \times
F^{\times},(1,\chi)}=0.$$ $\Gamma(\gd)$ has a finite number of $H$
orbits (it follows from Lemma \ref{EquivDef} and the introduction
of \cite{KR}). Hence it is enough to show that for any $x \in
\Gamma(\gd)$ we have
$$\Sc^*(Ad(H(F))x,Sym^k(CN_{Ad(H(F))x}^{\gd}))^{H(F) \times F^{\times},(1,\chi)}=0 \text{ for any } k.$$
Let $K:=\{(D_t(x),t^2) | t \in F^{\times}\} \subset (H(F)\times
F^{\times})_x.$

Note that
$$\Delta_{(H(F)\times
F^{\times})_x}((D_t(x),t^2))=|det(Ad(D_t(x))|_{\gd_x})| =
|t|^{\tr(ad(d(x))|_{\h_x})}.$$

By Lemma \ref{EigenInt} the eigenvalues of the action of
$(D_t(x),t^2)$ on $(Sym^k(\gd/[x,\h]))$ are of the form $t^l$
where $l$ is a non-positive integer.

Now by Frobenius reciprocity (Theorem \ref{Frob}) we have

\begin{multline*}
\Sc^*((H(F))x,Sym^k(CN_{Ad(H(F))x}^{\gd}))^{H(F) \times
F^{\times},(1,\chi)}=\\
=\Sc^*(\{x\},Sym^k(CN_{Ad(H(F))x,x}^{\gd})\otimes
\Delta_{H(F)\times F^{\times}}|_{(H(F)\times F^{\times})_x} \cdot
\Delta^{-1}_{(H(F)\times F^{\times})_x} \otimes
(1,\chi))^{(H(F)\times F^{\times})_x}=\\
=(Sym^k(\gd/[x,\h])\otimes \Delta_{(H(F)\times F^{\times})_x}
\otimes (1,\chi)^{-1} \otimes_{\R} \C )^{(H(F)\times F^{\times})_x}  \subset\\
\subset (Sym^k(\gd/[x,\h])\otimes \Delta_{(H(F)\times
F^{\times})_x} \otimes (1,\chi)^{-1}\otimes_{\R} \C)^{K}
\end{multline*}

which is zero since all the absolute values of the eigenvalues of
the action of any $(D_t(x),t^2) \in K$ on
$$Sym^k(\gd/[x,\h])\otimes \Delta_{(H(F)\times F^{\times})_x}
\otimes (1,\chi)^{-1}$$ are of the form $|t|^l$ where $l < 0$.
\end{proof}

\subsection{Regular symmetric pairs} \label{SecRegPairs}
$ $\\
In this subsection we will formulate a property which is weaker
than weakly linearly tame but still enables us to prove GK
property for good pairs.

\begin{definition}
Let $(G,H,\theta)$ be a symmetric pair. We call an element $g \in
G(F)$ \textbf{admissible} if\\
(i) $Ad(g)$ commutes with $\theta$ (or, equivalently, $s(g)\in Z(G)$) and \\
(ii) $Ad(g)|_{\g^{\sigma}}$ is $H$-admissible.
\end{definition}

\begin{definition} \label{DefReg}
We call a symmetric pair $(G,H,\theta)$ \textbf{regular} if for
any admissible $g \in G(F)$ such that
$\Sc^*(R(\g^{\sigma}))^{H(F)} \subset
\Sc^*(R(\g^{\sigma}))^{Ad(g)}$ we have
$$\Sc^*(Q(\gd))^{H(F)} \subset \Sc^*(Q(\gd))^{Ad(g)}.$$
\end{definition}

\begin{remark}
Clearly, every weakly linearly tame pair is regular.
\end{remark}

\begin{proposition}
A product of regular symmetric pairs is regular.
\end{proposition}
This is a direct corollary from Proposition \ref{Product}.

The goal of this subsection is to prove the following theorem.

\begin{theorem} \label{GoodHerRegGK}
Let $(G,H,\theta)$ be a good symmetric pair such that all its
descendants are regular. Then it is a GK pair.
\end{theorem}

We will need several definitions and lemmas.

\begin{definition}
Let $(G,H,\theta)$ be a symmetric pair. $g \in G$ is called
\textbf{normal} if $\sigma(g)g=g\sigma(g)$.
\end{definition}

The following lemma is straightforward.
\begin{lemma}
Let $(G,H,\theta)$ be a symmetric pair. Let $O \subset G(F)$ be an $H(F) \times H(F)$ orbit.\\
(i) If $\sigma(O)=O$ then there exists a normal element
$g \in O$. \\
(ii) Let $g \in G(F)$ be a normal element. Then there exists $h
\in H(F)$ such that $gh=hg=\sigma(g)$.
\end{lemma}

\begin{proof} $ $\\
(i). Let $g' \in O$. We know that $\sigma(g') = h_1 g' h_2$ where
$h_1,h_2 \in H(F)$. Let $g:=g' h_1$. Then
\begin{multline*}
\sigma(g)g = h_1^{-1} \sigma(g') g' h_1 = h_1^{-1} \sigma(g')
\sigma(\sigma(g')) h_1 = \\ = h_1^{-1} h_1 g' h_2
\sigma(h_1g'h_2))h_1= g' \sigma(g') = g' h_1 h_1^{-1} \sigma(g')=
g \sigma(g).
\end{multline*}
(ii) Follows from the fact that $g^{-1} \sigma(g) = \sigma(g)
g^{-1} \in H(F) $.
\end{proof}

\begin{notation}
Let $(G,H,\theta)$ be a symmetric pair. We denote
$\widetilde{H\times H}:= H\times H \rtimes \{1, \sigma\}$ where
$\sigma \cdot (h_1,h_2) = (\theta(h_2),\theta(h_1)) \cdot \sigma$.
The two-sided action of $H\times H$ on $G$ is extended to action
of $\widetilde{H\times H}$ in the natural way. We denote by $\chi$
the character of $\widetilde{H\times H}$ defined by
$\chi(\widetilde{H\times H} - H \times H)=\{-1\}$, $\chi(H \times
H)=\{1\}$.
\end{notation}

\begin{proposition}
Let $(G,H,\theta)$ be a good symmetric pair. Let $O \subset G(F)$ be a closed $H(F) \times H(F)$ orbit.\\
Then for any $g \in O$ there exist $\tau \in (\widetilde{H\times
H})_g(F) - (H\times H)_g(F)$ and $g' \in G_{s(g)}(F)$ such that
$Ad(g')$ commutes with $\theta$ on $G_{s(g)}$ and the action of
$\tau$ on $N_{O,g}^G$ corresponds via the isomorphism given by
Proposition \ref{PropDescend} to the adjoint action of $g'$ on
$\g_{s(g)}^{\sigma}$.
\end{proposition}
\begin{proof}
Clearly, if the statement holds for some $g\in O$ then it holds
for any $g\in O$.

Let $g \in O$ be a normal element. Let $h\in H(F)$ be such that
$gh=hg=\sigma(g)$. Let $\tau:= (h^{-1},1) \cdot \sigma$.
Evidently, $\tau\in (\widetilde{H\times H})_g(F) - (H\times
H)_g(F)$. Consider $d\tau_g: T_gG \to T_gG$. It corresponds via
the identification $dg: \g \cong T_gG$ to some $A:\g \to \g$.
Clearly, $A = da$ where $a:G \to G$ is defined by $a(\alpha) =
g^{-1}h^{-1}\sigma(g\alpha)$. However,
$g^{-1}h^{-1}\sigma(g\alpha) = \theta(g) \sigma(\alpha)
\theta(g)^{-1}.$ Hence $A = Ad(\theta(g)) \circ \sigma$. Let $B$
be a non-degenerate $G$-invariant $\sigma$-invariant symmetric
form on $\g$. By Theorem \ref{RegFun}, $A$ preserves $B$.
Therefore $\tau$ corresponds to $A|_{\g_{s(g)}^{\sigma}}$ via the
isomorphism given by Proposition \ref{PropDescend}. However,
$\sigma$ is trivial on $\g_{s(g)}^{\sigma}$ and hence
$A|_{\g_{s(g)}^{\sigma}}=Ad(\theta(g))|_{\g_{s(g)}^{\sigma}}$.
Since $g$ is normal, $\theta(g) \in G_{s(g)}$. It is easy to see
that $Ad(\theta(g))$ commutes with $\theta$ on $G_{s(g)}$. Hence
we take $g':=\theta(g)$.
\end{proof}

The last proposition implies Theorem \ref{GoodHerRegGK}. This
implication is proven in the same way as Theorem \ref{Invol_HC}.

\subsection{Conjectures} \label{conj}
\begin{conjecture}[van Dijk] \label{ConjDijk}
If $F=\C$, any connected symmetric pair is a Gelfand pair (GP3,
see Definition \ref{GPs} below).
\end{conjecture}
By theorem \ref{DistCrit} it follows from the following conjecture.
\begin{conjecture}\label{ConjDijkGK}
If $F=\C$, any connected symmetric pair is a GK pair.
\end{conjecture}
By Corollary \ref{ComplexGood} it follows from the following more
general conjecture.
\begin{conjecture} \label{ConjGKGood}
Every good symmetric pair is a GK pair.
\end{conjecture}
which in turn follows (by Theorem \ref{GoodHerRegGK})  from the
following one.
\begin{conjecture} \label{ConjAllReg}
Any symmetric pair is regular.
\end{conjecture}
An indirect evidence for this conjecture is that one can show that
every GK pair is regular.
\begin{remark}
It is well known that if $F$ is archimedean, $G$ is connected and
$H$ is compact then the pair $(G,H,\theta)$ is good, Gelfand (GP1,
see Definition \ref{GPs} below) and in fact also GK.
\end{remark}

\begin{remark}
In general, not every symmetric pair is good. For example,
$(SL_2(\R),T)$ where $T$ is the split torus. Also, it is not a
Gelfand pair (even not GP3, see Definition \ref{GPs} below).
\end{remark}

\begin{remark}
We do not believe that any symmetric pair is special. However, in
the next subsection we will prove that certain symmetric pairs are
special.
\end{remark}

\subsection{The pairs $(G \times G, \Delta G)$ and $(G_{E/F}, G)$
are tame} \label{Sec2RegPairs}

\begin{notation}
Let $E$ be a quadratic extension of $F$. Let $G$ be an algebraic
group defined over $F$. We denote by $G_{E/F}$ the canonical
algebraic group defined over $F$ such that $G_{E/F}(F)=G(E)$.
\end{notation}

In this section we will prove the following theorem.

\begin{theorem}\label{2RegPairs}
Let $G$ be a reductive group.\\
(i)Consider the involution $\theta$ of $G\times G$ given by
$\theta((g,h)):= (h,g)$. Its fixed points form the diagonal
subgroup $\Delta G$. Then the symmetric pair $(G \times G, \Delta
G, \theta)$ is
tame.\\
(ii) Let $E$ be a quadratic extension of $F$. Consider the
involution $\gamma$ of $G_{E/F}$ given by the nontrivial element
of $Gal(E/F)$. Its fixed points form $G$. Then the symmetric pair
$(G_{E/F}, G, \gamma)$ is tame.
\end{theorem}

\begin{corollary} \label{GCongTame}
Let $G$ be a reductive group. Then the adjoint action of $G$ on
itself is tame. In particular, every conjugation invariant
distribution on $GL_n(F)$ is transposition invariant \footnote{In
the non-archimedean case, the later is a classical result of
Gelfand and Kazhdan, see \cite{GK}.}.
\end{corollary}

For the proof of the theorem we will need the following
straightforward lemma.
\begin{lemma} $ $\\
(i) Every descendant of $(G \times G, \Delta G, \theta)$ is of the
form $(H \times H, \Delta H, \theta)$ for some reductive group
$H$.\\
\noindent (ii) Every descendant of $(G_{E/F}, G, \gamma)$ is of
the form $(H_{E/F}, H, \gamma)$ for some reductive group $H$.
\end{lemma}

Now Theorem \ref{2RegPairs} follows from the following theorem.

\begin{theorem}
The pairs $(G \times G, \Delta G, \theta)$ and  $(G_{E/F}, G,
\gamma)$ are special for any reductive group $G$.
\end{theorem}

By the speciality criterion (Proposition \ref{SpecCrit}) this
theorem follows from the following lemma.

\begin{lemma}
Let $\g$ be a semisimple Lie algebra. Let $\{e,h,f\} \subset \g$
be an $sl_2$ triple. Then $tr(Ad(h)|_{{\g}_e})$ is an integer
smaller than $dim\g$.
\end{lemma}
\begin{proof}
Consider $\g$ as a representation of $sl_2$ via the triple
$(e,h,f)$. Decompose it into irreducible representations
$\g=\bigoplus V_i$. Let $\lambda_i$ be the highest weights of
$V_i$. Clearly $$tr(Ad(h)|_{{\g}_e})= \sum \lambda_i \text{ and }
dim\g= \sum (\lambda_i+1).$$
\end{proof}

\section{Applications to Gelfand pairs} \label{Gel}
\subsection{Preliminaries on Gelfand pairs and distributional criteria}
$ $\\
In this section we recall a technique due to Gelfand and Kazhdan
which allows to deduce statements in representation theory from
statements on invariant distributions. For more detailed
description see \cite{AGS1}, section 2.

\begin{definition}
Let $G$ be a reductive group. By an \textbf{admissible
representation of} $G$ we mean an admissible representation of
$G(F)$ if $F$ is non-archimedean (see \cite{BZ}) and admissible
smooth \Fre representation of $G(F)$ if $F$ is archimedean.
\end{definition}

We now introduce three notions of Gelfand pair.

\begin{definition}\label{GPs}
Let $H \subset G$ be a pair of reductive groups.
\begin{itemize}
\item We say that $(G,H)$ satisfy {\bf GP1} if for any irreducible
admissible representation $(\pi,E)$ of $G$
we have
$$dim Hom_{H(F)}(E,\cc) \leq 1$$

\item We say that $(G,H)$ satisfy {\bf GP2} if for any irreducible
admissible representation $(\pi,E)$ of $G$
we have
$$dim Hom_{H(F)}(E,\cc) \cdot dim Hom_{H}(\widetilde{E},\cc)\leq
1$$

\item We say that $(G,H)$ satisfy {\bf GP3} if for any irreducible
{\bf unitary} representation $(\pi,\mathcal{H})$ of $G(F)$ on a
Hilbert space $\mathcal{H}$ we have
$$dim Hom_{H(F)}(\mathcal{H}^{\infty},\cc) \leq 1.$$
\end{itemize}

\end{definition}
Property GP1 was established by Gelfand and Kazhdan in certain
$p$-adic cases (see \cite{GK}). Property GP2 was introduced in
\cite{Gross} in the $p$-adic setting. Property GP3 was studied
extensively by various authors under the name {\bf generalized
Gelfand pair} both in the real and $p$-adic settings (see e.g.
\cite{vD-P}, \cite{vD}, \cite{Bos-vD}).

We have the following straightforward proposition.

\begin{proposition}
$GP1 \Rightarrow GP2 \Rightarrow GP3.$
\end{proposition}

We will use the following theorem from \cite{AGS1} which is a
version of a classical theorem of Gelfand and Kazhdan (see
\cite{GK}).

\begin{theorem}\label{DistCrit}
Let $H \subset G$ be reductive groups and let $\tau$ be an
involutive anti-automorphism of $G$ and assume that $\tau(H)=H$.
Suppose $\tau(\xi)=\xi$ for all bi $H(F)$-invariant Schwartz
distributions $\xi$ on $G(F)$. Then $(G,H)$ satisfies GP2.
\end{theorem}

In some cases, GP2 is equivalent to GP1. For example, see
corollary \ref{GKCor} below.

\subsection{Applications to Gelfand pairs}

\begin{theorem}\label{GKRep}
Let $G$ be reductive group and let $\sigma$ be an
$Ad(G)$-admissible anti-automorphism of $G$. Let $\theta$ be the
automorphism of $G$ defined by $\theta(g):=\sigma(g^{-1})$. Let
$(\pi, E)$ be an irreducible admissible representation of $G$.

Then $\widetilde{E} \cong E^{\theta}$, where $\widetilde{E}$
denotes the smooth contragredient representation and $E^{\theta}$
is $E$ twisted by $\theta$.

\end{theorem}
\begin{proof}
By Theorem 8.1.5 in \cite{Wal1}, it is enough to prove that the
characters of  $\widetilde{E}$ and $E^{\theta}$ are identical.
This follows from corollary \ref{GCongTame}.
\end{proof}

\begin{remark}
This theorem has an alternative proof using Harish-Chandra
regularity theorem, which says that character of an admissible
representation is a locally integrable function.
\end{remark}

\begin{corollary} \label{GKCor}
Let $H \subset G$ be reductive groups and let $\tau$ be an
$Ad(G)$-admissible anti-automorphism of $G$ such that $\tau(H)=H$.
Then $GP1$ is equivalent to $GP2$ for the pair $(G,H)$.
\end{corollary}

\begin{theorem}\label{Flicker}
Let $E$ be a quadratic extension of $F$. Then the pair $(GL_n(E),
GL_n(F))$ satisfies GP1.
\end{theorem}
For non-archimedean $F$ this theorem is proven in \cite{Fli}.
\begin{proof}
By theorem \ref{2RegPairs} this pair is tame. Hence it is enough
to show that this symmetric pair is good. This follows from the
fact that for any semisimple $x \in GL_n(E)^{\sigma}$ we have
$H^1(F,(GL_n)_x)=0$. Here we consider the adjoint action of $GL_n$
on itself.
\end{proof}

\appendix

\section{Localization principle}\label{SecRedLocPrin}
$ \quad \quad \quad \quad  \quad \quad$ by Avraham Aizenbud,
Dmitry Gourevitch
and Eitan Sayag\\\\
\setcounter{lemma}{0}
In this appendix we formulate and prove localization principle in
the case of a reductive group $G$ acting on a smooth affine
variety $X$. This is relevant only over archimedean $F$ since for
$l$-spaces, a more general version of this principle has been
proven in \cite{Ber}.

In \cite{AGS2}, we formulated localization principle in the
setting of differential geometry. Currently we do not have a proof
of this principle in such setting. Now we present a proof in the
case of a reductive group $G$ acting on a smooth affine variety
$X$. This generality is wide enough for all applications we had up
to now, including the one in \cite{AGS2}.


\begin{theorem}[Localization principle] \label{LocPrin2}
Let a reductive group $G$ act on a smooth algebraic variety $X$.
Let $Y$ be an algebraic variety and $\phi:X \to Y$ be an affine
algebraic $G$-invariant map. Let $\chi$ be a character of $G(F)$.
Suppose that for any $y \in Y(F)$ we have
$\cD_{X(F)}(\phi(F)^{-1}(y))^{G(F),\chi}=0$. Then
$\cD(X(F))^{G(F),\chi}=0$.
\end{theorem}

\begin{proof}
Clearly, it is enough to prove for the case when $X$ is affine, $Y
= X/G$ and $\phi = \pi_X(F)$. By the generalized Harish-Chandra
descent (Corollary \ref{Strong_HC_Cor}), it is enough to prove
that for any $G$-semisimple $x\in X(F)$, we have
$$\cD_{N_{Gx,x}^X(F)}(\Gamma(N_{Gx,x}^X))^{G_x(F),\chi} = 0.$$

Let $(U,p,\psi,S,N)$ be an analytic Luna slice at $x$. Clearly,
$$\cD_{N_{Gx,x}^X(F)}(\Gamma(N_{Gx,x}^X))^{G_x(F),\chi} \cong
\cD_{\psi(S)}(\Gamma(N_{Gx,x}^X))^{G_x(F),\chi}\cong
\cD_{S}(\psi^{-1}(\Gamma(N_{Gx,x}^X)))^{G_x(F),\chi}.$$ By
Frobenius reciprocity,
$$\cD_{S}(\psi^{-1}(\Gamma(N_{Gx,x}^X)))^{G_x(F),\chi}=\cD_{U}(G(F)\psi^{-1}(\Gamma(N_{Gx,x}^X)))^{G(F),\chi}.$$

By lemma \ref{Gamma}, $$G(F)\psi^{-1}(\Gamma(N_{Gx,x}^X))= \{y \in
X(F)| x \in \overline{G(F)y} \}.$$
Hence by Corollary \ref{EquivClassClosed},
$G(F)\psi^{-1}(\Gamma(N_{Gx,x}^X))$ is closed in $X(F)$. Hence
$$\cD_{U}(G(F)\psi^{-1}(\Gamma(N_{Gx,x}^X)))^{G(F),\chi}
 =
\cD_{X(F)}(G(F)\psi^{-1}(\Gamma(N_{Gx,x}^X)))^{G(F),\chi}.$$ Now,
$$G(F)\psi^{-1}(\Gamma(N_{Gx,x}^X)) \subset
\pi_X(F)^{-1}(\pi_X(F)(x))$$ and we are given
$$\cD_{X(F)}(\pi_X(F)^{-1}(\pi_X(F)(x)))^{G(F),\chi}=0$$ for any
$G$-semisimple $x$.
\end{proof}

\begin{remark} \label{LocPrinS}
An analogous statement holds for Schwartz distributions and the
proof is the same.
\end{remark}

\begin{corollary} \label{LocPrinSub}
Let a reductive group $G$ act on a smooth algebraic variety $X$.
Let $Y$ be an algebraic variety and $\phi:X \to Y$ be an affine
algebraic $G$-invariant submersion. Suppose that for any $y \in
Y(F)$ we have $\Sc^*(\phi^{-1}(y))^{G(F),\chi}=0$. Then
$\cD(X(F))^{G(F),\chi}=0$.
\end{corollary}

\begin{proof}
%
For any $y \in Y(F)$, denote $X(F)_y:=(\phi^{-1}(y))(F)$. Since
$\phi$ is a submersion, for any $y \in Y(F)$ the set $X(F)_y$ is a
smooth manifold. Moreover, $d\phi$ defines an isomorphism between
$N_{X(F)_y,z}^{X(F)}$ and $T_{Y(F),y}$ for any $z \in X(F)_y$.
Hence the bundle $CN_{X(F)_y}^{X(F)}$ is a trivial
$G(F)$-equivariant bundle.

We know that $$\Sc^*(X(F)_y)^{G(F),\chi}=0.$$ Therefore for any
$k$, we have
$$\Sc^*(X(F)_y,\Sym^k(CN_{X(F)_y}^{X(F)}))^{G(F),\chi}=0.$$ Thus by
Theorem \ref{NashFilt}, $\Sc^*_{X(F)}(X(F)_y)^{G(F),\chi}=0$. Now,
by Theorem \ref{LocPrin2} (and Remark \ref{LocPrinS}) this implies
that $\Sc^*(X(F))^{G(F),\chi}=0$. Finally, by Theorem
\ref{NoSNoDist} this implies $\cD(X(F))^{G(F),\chi}=0$.
\end{proof}

\begin{remark} \label{RemLocVectSys}
Theorem \ref{LocPrin} and Corollary \ref{LocPrinSub} have obvious
generalizations to constant vector systems, and the same proofs
hold.
\end{remark}

\section{Algebraic geometry over local fields} \label{AppLocField}
\subsection{Implicit function theorems}
\label{AppSub}

\begin{definition}
An analytic map $\phi:M \to N$ is called \textbf{\et map} if
$d_x\phi: T_xM \to T_xN$ is an isomorphism for any $x\in M$. An
analytic map $\phi:M \to N$ is called \textbf{submersion} if
$d_x\phi: T_xM \to T_xN$ is onto for any $x\in M$.
\end{definition}

We will use the following version of the inverse function theorem.

\begin{theorem} \label{InvFunct}
Let $\phi: M \to N$ be an \et map of analytic manifolds. Then it
is locally an isomorphism.
\end{theorem}
For proof see e.g. \cite{Ser}, Theorem 2 in section 9 of Chapter
III in part II.

\begin{corollary} \label{EtLocIs}
Let $\phi: X \to Y$ be a morphism of (not necessarily smooth)
algebraic varieties. Suppose that $\phi$ is \et at $x \in X(F)$.

Then there exists an open neighborhood $U \subset X(F)$ of $x$
such that $\phi|_U$ is a homeomorphism to its open image in
$Y(F)$.
\end{corollary}
For proof see e.g. \cite{Mum}, Chapter III, section 5, proof of
Corolary 2. There, the proof is given for the case $F=\C$ but it
works in the general case.

\begin{remark}
If $F$ is archimedean then one can choose $U$ to be
semi-algebraic.
\end{remark}

The following proposition is well known (see e.g. section 10 of
Chapter III in part II of \cite{Ser}).
\begin{proposition}
Any submersion $\phi:M \to N$ is open.
\end{proposition}

\begin{corollary}
Lemma \ref{OrbitIsOpen} holds. Namely, for any algebraic group $G$
and a closed algebraic subgroup $H \subset G$ the subset
$G(F)/H(F)$ is open and closed in $(G/H)(F)$.
\end{corollary}

\begin{proof}
Consider the map $\phi: G(F) \to (G/H)(F)$ defined by
$\phi(g)=gH$. Clearly, it is a submersion and its image is exactly
$G(F)/H(F)$. Hence, $G(F)/H(F)$ is open. Since there is a finite
number of $G(F)$ orbits in $(G/H)(F)$ and each of them is open for
the same reason, $G(F)/H(F)$ is also closed.
\end{proof}

\subsection{Luna slice theorem} \label{AppLun}
$ $\\
In this subsection we formulate Luna slice theorem and show how it
implies Theorem \ref{LocLuna}. For a survey on Luna slice theorem
we refer the reader to \cite{Dre} and the original paper
\cite{Lun}.

\begin{definition}
Let a reductive group $G$ act on affine varieties $X$ and $Y$. A
$G$-equivariant algebraic map $\phi:X \to Y$ is called
\textbf{strongly \et} if\\
(i) $\phi/G:X/G \to Y/G$ is \et\\
(ii) $\phi$ and the quotient morphism $\pi_X: X \to X/G$ induce a
$G$-isomorphism $X \cong Y \times _{Y/G}X/G$.
\end{definition}

\begin{definition}
Let $G$ be a reductive  group and $H$ be a closed reductive
subgroup. Suppose that $H$ acts on an affine variety $X$. Then $G
\times _{H} X$ denotes $(G\times X)/H$ with respect to the action
$h(g,x)=(gh^{-1},hx)$.
\end{definition}

\begin{theorem}[Luna slice theorem] \label{Luna}
Let a reductive group $G$ act on a smooth affine variety $X$. Let
$x \in X$ be $G$-semisimple.

Then there exists a locally closed smooth affine $G_x$-invariant
subvariety $Z \ni x$ of $X$ and a strongly \et algebraic map of
$G_x$ spaces $\nu: Z \to N_{Gx,x}^X$ such that the $G$-morphism
$\phi : G \times_{G_x} Z \to X$ induced by the action of $G$ on
$X$ is strongly \et.
\end{theorem}
\begin{proof}
It follows from Proposition 4.18, lemma 5.1 and theorems 5.2 and
5.3 in \cite{Dre}, noting that one can choose $Z$ and $\nu$ (in
our notations) to be defined over $F$.
\end{proof}

\begin{corollary}
Theorem \ref{LocLuna} holds. Namely:\\
Let a reductive group $G$ act on a smooth affine variety $X$. Let
$x \in X(F)$ be $G$-semisimple.

Then there exist\\
(i) an open $G(F)$-invariant B-analytic neighborhood $U$ of
$G(F)x$ in $X(F)$ with a
$G$-equivariant B-analytic retract $p:U \to G(F)x$ and\\
(ii) a $G_x$-equivariant B-analytic embedding $\psi:p^{-1}(x)
\hookrightarrow N_{Gx,x}^{X}(F)$ with open saturated image such
that $\psi(x)=0$.
\end{corollary}
\begin{proof}
Let $Z$, $\phi$ and $\nu$ be as in the last theorem.

Let $Z':= Z/G_x \cong (G \times_{G_x} Z)/G$ and $X':=X/G$.
Consider the natural map $\phi': Z'(F) \to X'(F)$. By Corollary
\ref{EtLocIs} there exists a neighborhood $S' \subset Z'(F)$ of
$\pi_Z(x)$ such that $\phi'|_{S'}$ is a homeomorphism to its open
image.

Consider the natural map $\nu': Z'(F) \to N_{Gx,x}^X/G_x(F)$. Let
$S'' \subset Z(F)$ be a neighborhood of $\pi_Z(x)$ such that
$\nu'|_{S''}$ is an isomorphism to its open image. In case that F
is archimedean we choose $S'$ and $S''$ to be semi-algebraic.

Let $S:=\pi_Z^{-1}(S''\cap S') \cap Z(F)$. Clearly, $S$ is
B-analytic.

Let $\rho: (G \times_{G_x}Z)(F) \to Z'(F)$ be the natural
projection. Let $O= \rho^{-1}(S'' \cap S')$. Let $q:O \to
G/G_x(F)$ be the natural projection. Let $O':=q^{-1}(G(F)/G_x(F))$
and $q':=q|_{O'}$.

Now put $U := \phi(O')$ and put $p:U \to G(F)x$ be the morphism
that corresponds to $q'$. Note that $p^{-1}(x) \cong S$ and put
$\psi:p^{-1}(x) \to N_{Gx,x}^X(F)$ to be the imbedding that
corresponds to $\nu|_S$.
\end{proof}

\section{Schwartz distributions on Nash manifolds} \label{AppSubFrob}

\subsection{Preliminaries and notations}
$ $\\
In this appendix we will prove some properties of $K$-equivariant
Schwartz distributions on Nash manifolds. We work in the notations
of \cite{AG1}, where one can read on Nash manifolds and Schwartz
distributions over them. More detailed references on Nash
manifolds are \cite{BCR} and \cite{Shi}.

Nash manifolds are equipped with \textbf{restricted topology}.
This is the topology in which open sets are open semi-algebraic
sets. This is not a topology in the classical sense of the word as
infinite unions of open sets are not necessary open sets in the
restricted topology. However, finite unions of open sets are open
sets and therefore in the restricted topology we consider only
finite covers. In particular, if $E \to M$ is a Nash vector bundle
it means that there exists a \underline{finite} open cover $U_i$
of $M$ such that $E|_{U_i}$ is trivial.

\begin{notation}
Let $M$ be a Nash manifold. We denote by $D_M$ the Nash bundle of
densities on $M$. It is the natural bundle whose smooth sections
are smooth measures, for precise definition see e.g. \cite{AG1}.
\end{notation}

An important property of Nash manifolds is
\begin{theorem}[Local triviality of Nash manifolds.] \label{loctriv}
Any Nash manifold can be covered by finite number
of open submanifolds Nash diffeomorphic to $\R^n$.
\end{theorem}
For proof see theorem I.5.12 in \cite{Shi}.

\begin{definition}
Let $M$ be a Nash manifold. We denote by $\G(M):= \Sc^*(M,D_M)$
the \textbf{space of Schwartz generalized functions} on $M$.
Similarly, for a Nash bundle $E \to M$ we denote by $\G(M,E):=
\Sc^*(M,E^* \otimes D_M)$ the \textbf{space of Schwartz
generalized sections} of $E$.

In the same way, for any smooth manifold $M$ we denote by
$C^{-\infty}(M):= \cD(M,D_M)$ the \textbf{space of generalized
functions} on $M$ and for a smooth bundle $E \to M$ we denote by
$C^{-\infty}(M,E):= \cD(M,E^* \otimes D_M)$ the \textbf{space of
generalized sections} of $E$.
\end{definition}

Usual $L^1$ functions can be interpreted as Schwartz generalized
functions but not as Schwartz distributions. We will need several
properties of Schwartz functions from \cite{AG1}.

\begin{property}\label{pClass}  $\Sc(\R ^n)$ = Classical
Schwartz functions on $\R ^n$.
\end{property}
For proof see theorem 4.1.3 in \cite{AG1}.

\begin{property}\label{pOpenSet}
Let $U \subset M$  be a (semi-algebraic) open subset, then
$$\Sc(U,E) \cong \{\phi \in \Sc(M,E)| \quad \phi \text{ is 0 on } M
\setminus U \text{ with all derivatives} \}.$$
\end{property}
For proof see theorem 5.4.3 in \cite{AG1}.

\begin{property}\label{pCosheaf}
Let $M$ be a Nash manifold. Let $M = \bigcup U_i$ be a finite open
cover of $M$. Then a function $f$ on $M$ is a Schwartz function if
and only if it can be written as $f= \sum \limits _{i=1}^n f_i$
where $f_i \in \Sc(U_i)$ (extended by zero to $M$).

Moreover, there exists a smooth partition of unity $1 =\sum
\limits _{i=1}^n \lambda_i$ such that for any Schwartz function $f
\in \Sc(M)$ the function $\lambda_i f$ is a Schwartz function on
$U_i$ (extended by zero to $M$).
\end{property}
For proof see section 5 in \cite{AG1}.

\begin{property}\label{pSheaf}
Let $M$ be a Nash manifold and $E$ be a Nash bundle over it. Let
$M = \bigcup U_i$ be a finite open cover of $M$. Let $\xi_i \in
\G(U_i,E)$ such that $\xi_i|_{U_j} = \xi_j|_{U_i}$. Then there
exists a unique $\xi \in \G(M,E)$ such that $\xi|_{U_i} = \xi_i$.
\end{property}
For proof see section 5 in \cite{AG1}.

We will also use the following notation.
\begin{notation}
Let $M$ be a metric space and $x \in M$. We denote by $B(x,r)$ the
open ball with center $x$ and radius $r$.
\end{notation}
\subsection{Submersion principle}

\begin{theorem} \label{SurSubSec}
Let $M$ and $N$ be Nash manifolds and $s:M \rightarrow N$ be a
surjective submersive Nash map. Then locally it has a Nash
section, i.e. there exists a finite open cover $N= \bigcup \limits
_{i=1}^k U_i$ such that $s$ has a Nash section on each $U_i$.
\end{theorem}
For proof see \cite{AG2}, theorem 2.4.16.

\begin{corollary} \label{EtLocIsNash}
An \et map $\phi:M \to N$ of Nash manifolds is locally an
isomorphism. That means that there exist a finite cover $M =
\bigcup U_i$ such that $\phi|_{U_i}$ is an isomorphism to its open
image.
\end{corollary}

\begin{theorem} \label{NashEquivSub}
Let $p:M \to N$ be a Nash submersion of Nash manifolds. Then there
exist a finite  open (semi-algebraic) cover $M = \bigcup U_i$ and
isomorphisms $\phi_i:U_i \cong W_i$ and $\psi_i:p(U_i) \cong V_i$
where $W_i\subset \R^{d_i}$ and $V_i \subset \R^{k_i}$ are open
(semi-algebraic) subsets, $k_i \leq d_i$ and $p|_{U_i}$ correspond
to the standard projections.
\end{theorem}
\begin{proof}
Without loss of generality we can assume that $N = \R^k$, $M$ is
an equidimensional closed submanifold of $\R^n$ of dimension $d$,
$d \geq k$, and $p$ is given by the standard projection $\R^n \to
\R^k$.

Let $\Omega$ be the set of all coordinate subspaces of $\R^n$ of
dimension $d$ which contain $N$. For any $V \in \Omega$ consider
the projection $pr:M \rightarrow V$. Define $U_V= \{x \in M | d_x
pr$ is an isomorphism $\}$. It is easy to see that $pr|_{U_V}$ is
\`{e}tale and $\{U_V\}_{V \in \Omega}$ gives a finite cover of
$M$. Now the theorem follows from the previous corollary
(Corollary \ref{EtLocIsNash}).
\end{proof}

\begin{theorem} \label{NashSub}
Let $\phi:M \to N$ be a Nash submersion of Nash manifolds. Let $E$
be a Nash bundle over $N$. Then\\
(i) there exists a unique continuous linear map
$\phi_*:\Sc(M,\phi^*(E)\otimes D_M) \to \Sc(N,E \otimes D_N)$ such
that for any $f \in \Sc(N,E^*)$ and $\mu \in \Sc(M,\phi^*(E)
\otimes D_M)$ we have $$\int_{x \in N} \langle f(x),\phi_*\mu(x)
\rangle = \int_{x \in M} \langle \phi^*f(x), \mu(x) \rangle.$$ In
particular, we mean that both integrals converge. \\
(ii) If $\phi$ is surjective then $\phi_*$ is surjective.
\end{theorem}
\begin{proof}
$ $\\ \indent (i)

Step 1. Proof for the case when $M= \R^n$, $N= \R^k$, $k \leq n$,
$\phi$ is the
standard projection and $E$ is trivial.\\
Fix Haar measure on $\R$ and identify $D_{\R^l}$ with the trivial
bundle for any $l$. Define $$\phi_*(f)(x):= \int _{y \in \R^{n-k}}
f(x,y)dy.$$ Convergence of the integral and the fact that
$\phi_*(f)$ is a Schwartz function follows from standard calculus.

Step 2. Proof for the case when  $M \subset \R^n$ and $N \subset
\R^k$ are open
(semi-algebraic) subsets, $\phi$ is the standard projection and $E$ is trivial.\\
Follows from the previous step and Property \ref{pOpenSet}.

Step 3. Proof for the case when  $E$ is trivial.\\
Follows from the previous step, Theorem \ref{NashEquivSub} and
partition of unity (Property \ref{pCosheaf}).

Step 4. Proof in the general case.\\
Follows from the previous step and partition of unity (Property
\ref{pCosheaf}).

(ii) The proof is the same as in (i) except of Step 2. Let us
prove (ii) in the case of Step 2. Again, fix Haar measure on $\R$
and identify $D_{\R^l}$ with the trivial bundle for any $l$. By
Theorem \ref{SurSubSec} and partition of unity (Property
\ref{pCosheaf}) we can assume that there exists a Nash section
$\nu:N \to M$. We can write $\nu$ in the form $\nu(x) = (x,s(x))$.

For any $x \in N$ define $R(x):= \sup\{r \in \R_{\geq 0}|
B(\nu(x),r) \subset M \}$. Clearly, $R$ is continuous and
positive. By Tarski - Seidenberg principle (see e.g. \cite{AG1},
theorem 2.2.3) it is semi-algebraic. Hence (by lemma A.2.1 in
\cite{AG1}) there exists a positive Nash function $r(x)$ such that
$r(x) < R(x)$. Let $\rho \in \Sc(\R^{n-k})$ such that $\rho$ is
supported in the unit ball and its integral is 1. Now let $f \in
\Sc(N)$. Let $g \in C^{\infty}(M)$ defined by $g(x,y):=
f(x)\rho((y-s(x))/r(x))/r(x)$ where $x \in N$ and $y \in
\R^{n-k}$. It is easy to see that $g \in \Sc(M)$ and $\phi_*g=f$.
\end{proof}

\begin{notation}
Let $\phi:M \to N$ be a Nash submersion of Nash manifolds. Let $E$
be a bundle on $N$. We denote by $\phi^*:\G(N,E) \to
\G(M,\phi^*(E))$ the dual map to $\phi_*$.
\end{notation}

\begin{remark}
Clearly, the map $\phi^*:\G(N,E) \to \G(M,\phi^*(E))$ extends to
the map $\phi^*:C^{-\infty}(N,E) \to C^{-\infty}(M,\phi^*(E))$
described in \cite{AGS1}, theorem A.0.4.
\end{remark}

\begin{proposition} \label{EnoughPull}
Let $\phi:M \to N$ be a surjective Nash submersion of Nash
manifolds. Let $E$ be a bundle on $N$. Let $\xi \in
C^{-\infty}(N)$. Suppose that $\phi^*(\xi) \in \G(M)$. Then $\xi
\in \G(N)$.
\end{proposition}
\begin{proof}
It follows from Theorem \ref{NashSub} and Banach open map theorem
(see theorem 2.11 in \cite{Rud}).
\end{proof}

\subsection{Frobenius reciprocity}
$ $\\
In this subsection we prove Frobenius reciprocity for Schwartz
functions on Nash manifolds.

\begin{proposition}
Let $M$ be a Nash manifold. Let $K$ be a Nash group. Let $E \to M$
be a Nash bundle. Consider the standard projection $p:K \times M
\to M$. Then the map $p^*:\G(M,E) \to \G(M \times K, p^*E)^K$ is
an isomorphism.
\end{proposition}
This proposition follows from Proposition 4.0.11 in \cite{AG2}.

\begin{corollary}
Let a Nash group $K$ act on a Nash manifold $M$. Let $E$ be a
$K$-equivariant Nash bundle over $M$. Let $N \subset M$ be a Nash
submanifold such that the action map $K \times N \to M$ is
submersive. Then there exists a canonical map $$HC: \G(M,E)^K \to
\G(N,E|_N).$$
\end{corollary}

\begin{theorem}
Let a Nash group $K$ act on a Nash manifold $M$. Let $N$ be a
$K$-transitive Nash manifold. Let $\phi:M \to N$ be a Nash
$K$-equivariant map.

Let $z \in N$ be a point and $M_z:= \phi^{-1}(z)$ be its fiber.
Let $K_z$ be the stabilizer of $z$ in $K$. Let $E$ be a
$K$-equivariant Nash vector bundle over $M$.

Then there exists a canonical isomorphism $$Fr:
\G(M_z,E|_{M_z})^{K_z} \cong \G(M,\E)^K.$$
\end{theorem}
\begin{proof}

Consider the map $a_z:K \to N$ given by $a_z(g)=gz$. It is a
submersion. Hence by Theorem \ref{SurSubSec} there exists a finite
open cover $N= \bigcup \limits _{i=1}^k U_i$ such that $a_z$ has a
Nash section $s_i$ on each $U_i$. This gives an isomorphism
$\phi^{-1}(U_i) \cong U_i \times M_z$ which defines a projection
$p: \phi^{-1}(U_i) \to M_z$. Let $\xi \in \G(M_z,E|_{M_z})^{K_z}$.
Denote $\xi_i:=p^* \xi$. Clearly it does not depend on the section
$s_i$. Hence $\xi_i|_{U_i \cap U_j}=\xi_j|_{U_i \cap U_j}$ and
hence by Property \ref{pSheaf} there exists $\eta \in \G(M,\E)$
such that $\eta|_{U_i} = \xi_i$. Clearly $\eta$ does not depend on
the choices. Hence we can define $Fr(\xi) = \eta$.

It is easy to see that the map $HC: \G(M,E)^K \to
\G(M_{z},E|_{M_{z}})$ described in the last corollary gives the
inverse map.
\end{proof}

Since our construction coincides with the construction of
Frobenius reciprocity for smooth manifolds (see e.g. \cite{AGS1},
theorem A.0.3) we obtain the following corollary.

\begin{corollary}
Part (ii) of Theorem \ref{Frob} holds.
\end{corollary}

\subsection{$K$-invariant distributions compactly supported modulo
$K$.} \label{KInvAreSchwartz}
$ $\\

In this subsection we prove Theorem \ref{CompSchwartz}. Let us
first remind its formulation.

\begin{theorem} \label{CompSchwartz2}
Let a Nash group $K$ act on a Nash manifold $M$. Let $E$ be a
$K$-equivariant Nash bundle over $M$. Let $\xi \in \cD(M,E)^K$
such that $\Supp(\xi)$ is Nashly compact modulo $K$. Then $\xi \in
\Sc^*(M,E)^K$.
\end{theorem}

For the proof we will need the following lemmas.

\begin{lemma} \label{RelComp}
Let $M$ be a Nash manifold. Let $C \subset M$ be a compact subset.
Then there exists a relatively compact open (semi-algebraic)
subset $U \subset M$ that includes $C$.
\end{lemma}

\begin{proof}
For any point $x \in C$ choose an affine chart, and let $U_x$ be
an open ball with center at $x$ inside this chart. Those $U_x$
give an open cover of $C$. Choose a finite subcover
$\{U_i\}_{i=1}^n$ and let $U:= \bigcup_{i=1}^n U_i$.
\end{proof}

\begin{lemma} \label{takayata}
Let $M$ be a Nash manifold. Let $E$ be a
Nash bundle over $M$. Let $U \subset M$ be a relatively
compact open (semi-algebraic) subset. Let $\xi \in \cD(M,E)$. Then
$\xi|_U \in \Sc^*(U,E|_U)$.
\end{lemma}

\begin{proof}
It follows from the fact that extension by zero $ext:\Sc(U,E|_U) \to
C_c^{\infty}(M,E)$ is a continuous map.
\end{proof}

\begin{proof}[Proof of Theorem \ref{CompSchwartz2}]
Let $Z \subset M$ be a semi-algebraic closed subset and $C \subset M$ be a compact subset such that $Supp(\xi) \subset Z \subset
KC$.

Let $U \supset C$ be as in Lemma \ref{RelComp}. Let $\xi' :=
\xi|_{KU}$. Since $\xi|_{M - Z}=0$, it is enough to show that $\xi'$ is Schwartz.

Consider the surjective submersion $m_U:K \times U \to
KU$. Let $$\xi'':=m_U^* (\xi') \in \cD(K \times U,m_U^*(E))^K.$$
By Proposition \ref{EnoughPull}, it is enough to show that $$\xi''
\in \Sc^*(K \times U,m_U^*(E)).$$
By Frobenius reciprocity, $\xi''$ corresponds to $\eta \in
\cD(U,E)$. It is enough to prove that $\eta \in \Sc^*(U,E)$.
Consider the submersion $m:K \times M \to M$ and let $$\xi''' :=
m^*(\xi) \in \cD(K \times M,m^*(E)).$$ By Frobenius reciprocity,
$\xi'''$ corresponds to $\eta' \in \cD(M,E)$. Clearly $\eta =
\eta'|_U$. Hence by Lemma \ref{takayata}, $\eta \in \Sc^*(U,E)$.
\end{proof}
\section{Proof of archimedean homogeneity theorem} \label{AppRealHom}
The goal of this appendix is to prove Theorem \ref{ArchHom} for
archimedean $F$. First we remind its formulation.

\setcounter{lemma}{0}

\begin{theorem} [archimedean homogeneity]
Let $V$ be a vector space over $F$. Let $B$ be a non-degenerate
symmetric bilinear form on $V$. Let $M$ be a Nash manifold. Let $L
\subset \Sc^*_{V\times M}(Z(B)\times M)$ be a non-zero subspace
such that $\forall \xi \in L $ we have $\Fou_B(\xi) \in L$ and $B
\xi \in L$ (here $B$ is interpreted as a quadratic form).

Then there exists a non-zero distribution $\xi \in L$ which is
adapted to $B$.
\end{theorem}
\noindent Till the end of the section we assume that $F$ is
archimedean and we fix $V$ and $B$.

First we will need some facts about the Weil representation. For a
survey on the Weil representation in the archimedean case we refer
the reader to \cite{RS1}, section 1.

\begin{enumerate}
\item There exists a unique (infinitesimal) action
$\pi$ of $sl_2(F)$ on $\Sc^*(V)$ such that\\
(i) $\pi(\begin{pmatrix}
  0 & 1 \\
  0 & 0
\end{pmatrix}) \xi = -i \pi Re(B)
\xi$ and $\pi(\begin{pmatrix}
  0 & 0 \\
  -1 & 0
\end{pmatrix}) \xi = -\Fou_B^{-1} (i \pi Re(B)
\Fou_B(\xi))$.\\
(ii) If $F=\C$ then $\pi(\begin{pmatrix}
  0 & i \\
  0 & 0
\end{pmatrix})  = \pi(\begin{pmatrix}
  0 & 0 \\
  -i & 0
\end{pmatrix})=0$

\item It can be lifted to an action of the metaplectic group
$Mp(2,F)$.

We will denote this action by $\Pi$.

\item In case $F=\C$ we have $Mp(2,F)=SL_2(F)$ and in case $F=\R$ the
group $Mp(2,F)$ is a connected 2-folded covering of $SL_2(F)$. We
will denote by $\eps \in Mp(2,F)$ the element of order 2 that
satisfies $SL_2(F)=Mp(2,F)/\{1, \eps\}.$

\item In case $F=\R$ we have $\Pi(\eps)=(-1)^{dim V}$ and therefore if $dim
V$ is even then $\Pi$ factors through $SL_2(F)$ and if $dim V$ is
odd then no nontrivial subrepresentation of $\Pi$ factors through
$SL_2(F)$. In particular if $dim V$ is odd then $\Pi$ has no
nontrivial finite dimensional representations, since every finite
dimensional representation of $sl_2$ has a unique lifting both to
$SL_2(F)$ and to $Mp(2,F)$. \label{FinDimSubrep}

\item In case $F=\C$ or in case $dim V$ is even we have $\Pi(\begin{pmatrix}
  t & 0 \\
  0 & t^{-1}
\end{pmatrix}) \xi=\delta^{-1}(t)|t|^{-dimV/2} \rho(t) \xi$ and $\Pi(\begin{pmatrix}
  0 & 1 \\
  -1 & 0
\end{pmatrix}) \xi=
\gamma(B)^{-1} \Fou_B \xi.$
\end{enumerate}

We also need the following straightforward lemma.

\begin{lemma}
Let $(\Lambda, L)$ be a continuous finite dimensional
representation of $SL_2(\R)$.  Then there exists a non-zero $\xi
\in L$ such that either $$ \Lambda(\begin{pmatrix}
  t & 0 \\
  0 & t^{-1}
\end{pmatrix})\xi=\xi \text{ and } \Lambda(\begin{pmatrix}
  0 & 1 \\
  -1 & 0
\end{pmatrix})\xi \text{is proportional to } \xi$$  or
$$\Lambda(\begin{pmatrix}
  t & 0 \\
  0 & t^{-1}
\end{pmatrix})\xi=t \xi.$$
\end{lemma}

Now we are ready to prove the theorem.
\begin{proof}[Proof of Theorem \ref{ArchHom}]
Without loss of generality assume $M=pt$.

Let $\xi \in L$ be a non-zero distribution. Let $L':=
U_{\C}(sl_2(\R)) \xi \subset L$. Here, $U_{\C}$ means the
complexified universal enveloping algebra.

It is easy to see that $L' \subset \Sc^*(V)$ is finite dimensional
(see Lemma \ref{FinDim} below). Clearly, $L'$ is also
$sl_2(F)$-invariant and hence is also a subrepresentation of
$\Pi$. Therefore by Fact (\ref{FinDimSubrep}), $F=\C$ or $dim V$
is even. Hence $\Pi$ factors through $SL_2(F)$.

Now by the lemma there exists $\xi' \in L'$ which is $B$ adapted.
\end{proof}

\begin{lemma}\label{FinDim}
Let $V$ be a representation of $sl_2$. Let $v \in V$ be a vector
such that $e^k v=f^n v=0$ for some $n,k$. Then the representation
generated by $v$ is finite dimensional.\footnote{For our purposes
it is enough to prove this lemma for k=1.}
\end{lemma}

\begin{proof}
The proof is by induction on k.\\\\
Base k=1:\\
It is easy to see that $$e^l f^l v=l!(\prod_{i=0}^{l-1}(h-i))v$$
for all $l$. It can be checked by a direct computation, and also
follows from the fact that $e^l f^l$ is of weight $0$, hence it
acts on the singular vector $v$ by its Harish Chandra projection
which is $HC(e^l f^l)=l! \prod_{i=0}^{l-1}(h-i)$.

Therefore $(\prod_{i=0}^{n-1}(h-i))v=0$.

Hence $W:=U_{\C}(h) v$ is finite dimensional and $h$ acts on it
semi-simply. Let $\{v_i\}_{i=1}^m$ be an eigenbasis of $h$ in $W$.
It is enough to show that $U_{\C}(sl_2)v_i$ is finite dimensional
for any $i$. Note that $e|_W=f^n|_W=0$. Now, $U_{\C}(sl_2)v_i$ is
finite dimensional by Poincare-Birkhoff-Witt
theorem.\\\\
Induction step:\\
Let $w:=e^{k-1}v$. Let us show that $f^{n+k-1} w=0$. Consider the
element $f^{n+k-1}e^{k-1}\in U_{\C}(sl_2)$. It is of weight $-2n$,
hence by Poincare-Birkhoff-Witt theorem it can be rewritten as a
combination of elements of the form $e^a h^b f^c$ such that $c-a=
n$ and hence $c \geq n$. Therefore $f^{n+k-1}e^{k-1} v=0$.

Now let $V_1:=U_{\C}(sl_2) v$ and $V_2:=U_{\C}(sl_2) w$. By the
base of the induction $V_2$ is finite dimensional, by the
induction hypotheses  $V_1/V_2$ is finite dimensional, hence $V_1$
is finite dimensional.
\end{proof}

\newpage

\section{Diagram} \label{Diag}

The following diagram illustrates the interrelations of various
properties of a symmetric pair $(G,H)$. On the non-trivial implications we put the numbers of the statements that prove them. \\\\
\small
\newlength{\gnat}

\setlength{\gnat}{10pt}
$  \xymatrix{& & & &
\framebox{\parbox{65pt}{For any \\ nilpotent $x \in \gd$
$\tr(ad(d(x))|_{\h_x})$ $< dim \g^{\sigma} $}}\ar@{=>}[d]^{\ref{SpecCrit}} & &\\
& & & &\framebox{\parbox{25pt}{special}}\ar@{=>}[d]^{\ref{SpecWeakReg}}& & \\
\framebox{\parbox{27pt}{regular}}& & &
&\framebox{\parbox{33pt}{weakly linearly
tame}}\ar@{=>}[llll]&\framebox{\parbox{53pt}{All the\\ descendants
are weakly \\ linearly tame}}
\ar@{=>}[dl]_{\ref{LinDes}}\ar@{=>}[ddl]_{\ref{LinDes}}\\
& &\framebox{\parbox{47pt}{For any\\ descendant $(G',H')$:\\
$H^1(F,H')$ is trivial }}
\ar@{=>}[d] & &\framebox{\parbox{33pt}{linearly tame}}\ar@{=>}[u]&\\
 \framebox{\parbox{49pt}{All the\\ descendants are regular}}\ar@{=}[r]&AND\ar@{=}[r]\ar@{=>}
 [rd]_{\ref{GoodHerRegGK}}&\framebox{\parbox{20pt}{good}}\ar@{=}[r]&AND\ar@{=}[r]\ar@{=>}[ld] &\framebox{\parbox{33pt}{tame}}
 &\\
&&\framebox{\parbox{20pt}{GK}}\ar@/_5pc/@{=>}[dd]_{\ref{DistCrit}}&&&&\\
&&\framebox{\parbox{20pt}{GP3}}&&&\\
&&\framebox{\parbox{20pt}{GP2}}\ar@{=>}[u]\ar@{=}[r]&AND\ar@{=>}[ld]_{\ref{GKCor}}\ar@{=}[r]&\framebox{\parbox{55pt}
{$G$ has an\\ $Ad(G)$-admissible anti-automorphism that\\ preserves
$H$}}& \ar@{=>}[l] \framebox{\parbox{40pt}
{$G=GL_n$ and\\ $H=H^t$ }} \\
&&\framebox{\parbox{20pt}{GP1}}\ar@{=>}[u]&&&\\
 }$

\end{document}